\documentclass[11pt]{article}        
\usepackage{graphicx}
\begin{document}
\newtheorem{defn0}{Definition}[section]
\newtheorem{prop0}[defn0]{Proposition}
\newtheorem{thm0}[defn0]{Theorem}
\newtheorem{lemma0}[defn0]{Lemma}
\newtheorem{coro0}[defn0]{Corollary}
\newtheorem{exa}[defn0]{Example}
\def\rig#1{\smash{ \mathop{\longrightarrow}
    \limits^{#1}}}
\def\swar#1{\swarrow
   \rlap{$\vcenter{\hbox{$\scriptstyle#1$}}$}}
\def\lswar#1{\swarrow
   \llap{$\vcenter{\hbox{$\scriptstyle#1$}}$}}
\def\sear#1{\searrow
   \rlap{$\vcenter{\hbox{$\scriptstyle#1$}}$}}
\def\lsear#1{\searrow
   \llap{$\vcenter{\hbox{$\scriptstyle#1$}}$}}   
\def\near#1{\nearrow
   \rlap{$\vcenter{\hbox{$\scriptstyle#1$}}$}}   
\def\dow#1{\Big\downarrow
   \rlap{$\vcenter{\hbox{$\scriptstyle#1$}}$}}
\def\ldow#1{\Big\downarrow
   \llap{$\vcenter{\hbox{$\scriptstyle#1$}}$}}   
\def\up#1{\Big\uparrow
   \rlap{$\vcenter{\hbox{$\scriptstyle#1$}}$}}
\def\lef#1{\smash{ \mathop{\longleftarrow}
    \limits^{#1}}}
\newcommand{\defref}[1]{Def.~\ref{#1}}
\newcommand{\propref}[1]{Prop.~\ref{#1}}
\newcommand{\thmref}[1]{Thm.~\ref{#1}}
\newcommand{\lemref}[1]{Lemma~\ref{#1}}
\newcommand{\corref}[1]{Cor.~\ref{#1}}
\newcommand{\exref}[1]{Example~\ref{#1}}
\newcommand{\secref}[1]{Section~\ref{#1}}

\newcommand{\qedd}{\hfill\framebox[2mm]{\ }}

\author{Giorgio Ottaviani\ \ \ \ \ \ \ \ Elena Rubei}

\title{Quivers and the cohomology of\\ homogeneous vector bundles}
\date{}
\maketitle

\def\thefootnote{}
\footnotetext{ \hspace*{-0.36cm}MSC: Primary 14F05 Secondary 14M17, 32M15, 16G20 }

\begin{abstract}
We describe the cohomology groups of a homogeneous vector bundle
$E$ on any Hermitian symmetric variety $X=G/P$ of ADE type
as the cohomology
of a complex explicitly described. The main tool is the equivalence between
 the category of
 homogeneous bundles and the category of representations of a certain
 quiver ${\cal Q}_X$ with relations, whose vertices are the dominant 
weights of the reductive part of $P$. This equivalence was found in some
 cases by
Bondal, Kapranov and Hille and we find the appropriate relations
on any  Hermitian
 symmetric variety. 
\end{abstract}

\section{Introduction}
The Borel-Weil-Bott theorem computes the cohomology groups
of an irreducible homogeneous bundle on a rational homogeneous variety $X$. In this paper we
compute the cohomology groups of any homogeneous bundle (including the reducible ones)
on a symmetric Hermitian variety of ADE type. This class of varieties includes
grassmannians, quadrics of even dimension, spinor varieties, two exceptional cases and
products among all of them.

In order to compute the cohomology groups 
(see \thmref{filt1} ) we have to describe 
the homogeneous bundles as representations of a certain quiver ${\cal Q}_X$. 
The moduli spaces of such representations give
 moduli spaces of homogeneous bundles, that are introduced in \S 7 and seem to
 have an intrinsic interest.
 
 We describe now with some detail the background of this paper. 

Let $X=G/P$ be a rational homogeneous variety. It is known that the category
 of 
$G$-homogeneous bundles on $X$ 
is equivalent to the category $P$-$mod$ of
representations of $P$, and also to the category
${\cal P}$-$mod$ where ${\cal P}=Lie P$ (see for example \cite{B-K}). 
Since $P$ is not reductive, its representations are difficult to describe.
 In fact if $E$
 is a homogeneous bundle, it has a
 filtration
$0\subset E_1\subset\ldots\subset E_k=E$ where $E_i/E_{i-1}$ is 
irreducible, but the filtration does not split in general.

 Let $P=R\cdot N$ be the Levi decomposition, where $R$ is reductive and
 $N$ is nilpotent. At the level of Lie algebras this amounts to
 ${\cal P}={\cal R}\oplus
{\cal N}$ as vector spaces. Considering $E$ as $R$-module (and hence as ${\cal R}$
-module) we get the graded bundle $gr E=\oplus_i E_i/E_{i-1}$. 
The nilpotent radical
${\cal N}$ is an ${\cal R}$-module itself, with the adjoint action,
 corresponding to the bundle $\Omega^1_X$. The action of ${\cal P}$
over $E$ induces a $G$-equivariant map
$$\hspace*{1.5cm}\theta\colon\Omega^1_X\otimes gr E\rig{} gr E \qquad\qquad (*)$$
Our first result is that, when $X$ is a Hermitian symmetric variety, 
a morphism of ${\cal R}$-modules 
$\theta\colon\Omega^1_X\otimes F\rig{} F$ 
is induced by a ${\cal P}$-action 
if and only if $\theta\wedge\theta=0$ (see \thmref{starting}).

In analogy with \cite{Simp}, we call a completely reducible bundle 
$F$ 
endowed with such $\theta$ satisfying $\theta\wedge\theta=0$ a
 (homogeneous)
 Higgs bundle. So the category of $G$-homogeneous bundles turns out to be
 equivalent to the category
of Higgs bundles. In the pair $(F, \theta)$, $F$ encodes the discrete part 
and $\theta$ encodes the continuous part.

 By using Bott theorem we can prove that $ Hom(gr E\otimes\Omega^1_X, gr E)^G$ is isomorphic
 to $Ext^1(gr E, gr E)^G$ (see \thmref{homext}).
In this setting a reformulation of \thmref{starting}
implies that 
the set of ${\cal P}$-modules $E$ such that $gr E=F$
is in natural bijection with the set of $e\in Ext^1(F, F)^G$ such that $m(e)=0$
where $m$ is the quadratic Yoneda morphism
$Ext^1(F, F)^G\rig{}Ext^2(F, F)^G$.

Bondal and Kapranov had the remarkable idea that quivers are the appropriate
 tool to manage $P$-modules, indeed we state our results in the framework
 of quivers. 

A quiver ${\cal Q}_X$ is associated to any rational homogeneous variety
$X$. The points of ${\cal Q}_X$ are the 
dominant weights of $R$ and the arrows correspond to the weights of
 ${\cal N}$ in the action (*). Bondal, Kapranov
\cite{B-K} and Hille \cite{Hi2} proved that the category
of $G$-homogeneous bundles on $X$ is equivalent to the category of
 representations of ${\cal Q}_X$ with certain relations
 to be determined (see also \cite{King}). Hille in \cite{Hi2} proved that the relations in ${\cal Q}_X$ 
are quadratic if $X$ is Hermitian symmetric and found that the relations of the quiver constructed in 
\cite{B-K},
although essentially corrected,
were not properly stated in the case
of the Grassmannian of lines in ${\bf P}^3$
(see \exref{g13}). Then Hille showed that 
in ${\cal Q}_{{\bf P}^2}$ the relations correspond to the commutativity
of all square diagrams.  If $X$ is Hermitian symmetric we see that
the relations are consequences of
the condition $\theta\wedge\theta=0$.  This allows one to extend Hille's result
to ${\cal Q}_{{\bf P}^n}$ (see \corref{commutpn}).

\medskip
 
The second part of the paper is devoted to the computation of the cohomology.
The Borel-Weil-Bott theorem computes the cohomology groups of an
 irreducible bundle $E$ on $X$. In particular it says that $H^*(E)$
is an irreducible $G$-module.
It follows that for any $G$-homogeneous bundle $E$ there is a spectral
 sequence constructed by the filtration $gr E$  
abutting to the cohomology groups of $E$.
The main problem is that the maps occurring in the spectral sequence,
although they are equivariant, are difficult to control.
In fact most of the main open problems about rational homogeneous varieties,
 like the computation of syzygies of their projective embeddings,
 reduce to the computation of cohomology groups of certain homogeneous 
bundles (see the recent book \cite{We}).

Assume now that $X$ is  Hermitian symmetric of ADE type.
Thanks to the Borel-Weil-Bott theorem, and to the results of Kostant in
\cite{Ko}, we can divide the points
of ${\cal Q}_X$ into several chambers, separated by the hyperplanes
 containing the singular weights, that we call Bott
chambers. 
We consider the segments
connecting any point of ${\cal Q}_X$
with its mirror images in the adjacent Bott chambers,
and we define certain linear maps
 $c_i\colon H^i(gr E)\to H^{i+1}(gr E)$,
by composing the maps associated to the representation of ${\cal Q}_X$
corresponding to $E$, along these segments. 
We get a sequence
\[\ldots\rig{}H^i(gr E)\rig{c_i}H^{i+1}(gr E)\rig{c_{i+1}}\ldots\]

{\it 
Our main result (\thmref{filt1}) is that  
this sequence is a complex
 and its cohomology 
(as $G$-module) 
is the usual cohomology $H^i(X,E)$.}

The proof of this result is obtained by comparing the maps $c_i$ with the boundary maps.
In the case of projective spaces
the computation of $c_i$ can be done quite easily.
It is worth remarking that
the derived category of homogeneous bundles was described by Kapranov
in the last section of \cite{Ka}.  The quivers allow one
to refine that approach.

It turns out from our proof that the cohomology modules $H^i(E)$
are equipped with a natural filtration
\[0\subset H^i[1](E)\subset H^i[2](E)\subset\ldots\subset H^i[N](E)=H^i(E)\] 

\medskip

The last part of the paper deals with moduli spaces.
There is a 
notion of semistability of representations of quivers introduced 
in \cite{King} (see also \cite{Migl}) which is suitable to construct moduli
 spaces according to Mumford GIT. This notion of semistability turns
 out to be equivalent to the Mumford-Takemoto
semistability of the bundle and we get moduli spaces of
 $G$-homogeneous semistable bundles with fixed $gr E$.
More precisely, the choice of an ${\cal R}$-module
$F$ is equivalent to the choice of a dimension vector $\alpha$ as in
\cite{King}. All semistable $P$-modules $E$ such that $gr E=F$ are parametrized
 by a projective moduli space $M_X(\alpha)$.
The properties of such moduli spaces  probably deserve further study.

\smallskip

Finally we want to mention that some applications of this approach to the case of homogeneous bundles on ${\bf P}^2$
appear in \cite{O-R}.

\medskip
We sketch now the content of the sections.
In \S 3 we describe the equivalence of categories between
 $G$-homogeneous bundles and Higgs bundles. In \S 4 we recall the Borel-Weil-Bott theorem,
 in the form found by Kostant (\cite{Ko}), which is suitable for our purposes.
In \S 5 we construct in detail the quiver ${\cal Q}_X$ with its relations and we prove the equivalence
 between the category of homogeneous bundles and the category
 of representations of ${\cal Q}_X$.
In \S 6 we prove our main result about the cohomology groups. 
In \S 7 we consider
the moduli spaces $M_X(\alpha)$ and we compare some different
notions of stability.
In \S 8 we make explicit for Grassmannians the relations stated in \S 5
by using the Olver maps.
 
We thank Laurent Manivel 
 for several helpful comments, and Piotr Pragacz, who adviced us
to consider the Olver maps in order to make the relations explicit.

 \section{Notation and Preliminaries}
In all the paper
  let $G$ be a semisimple complex Lie group. We fix a Cartan
 subalgebra ${\cal H}$ in $Lie G$. 
Let $\Delta =\{\alpha_1,\ldots ,\alpha_n\}$ be a {\it fundamental system of 
simple roots}
for $Lie G$. A positive root is a linear combination
with nonnegative integral coefficients of the  simple roots.
The Killing product allows one to identify ${\cal H}$
with ${\cal H}^{\vee}$ and thus to define the Killing product
also on ${\cal H}^{\vee}$.
Let  $\{\lambda_1,\ldots ,\lambda_n\}$ be  the
{\it fundamental 
weights} corresponding to $\{\alpha_1,\ldots ,\alpha_n\}$,
 i.e. the elements of ${\cal H}^{\vee}$ 
such that $\frac{2(\lambda_i,\alpha_j)}
{(\alpha_j,\alpha_j)}=\delta_{ij}$
where $(\ ,\ )$ is the Killing product. 
Let $Z$ be the lattice generated by the fundamental weights.
The elements in $Z$ that are linear combination with nonnegative
coefficients of the fundamental weights are called the
dominant weights for $G$, and they are the maximal weights
 of the irreducible representations of $Lie G$.  In the ADE case
all roots have length $\sqrt{2}$.

 For any $W$ representation of $G$ we denote
 by $W^G$ its 
invariant part, that is the subspace of $W$ where $G$ acts trivially. 
If $V$ is an irreducible representation, we denote
$W^V:=Hom(V, W)^G\otimes V$.

If $\lambda\in Z$ we denote by $V_{\lambda}$ the irreducible
representation of $G$ with highest weight $\lambda$. In the case
$G=SL(n+1)$ to any $\lambda$ is associated a Young diagram.
Precisely if we have $\lambda=\sum_{i=1}^n{n_i}\lambda_i$, then we set $a_i=\sum_{j\ge i}n_j$
and we get the Young diagram with $a_i$ boxes in the 
$i$-th row. We use the notation where the first row is the top row.
The $n$-uple $a=(a_1,\ldots, a_n)$ is a partition
of $\sum a_i$ and it is customary to denote $V_{\lambda}$
as ${\cal S}^aV$. In particular ${\cal S}^2V=Sym^2 V$
and ${\cal S}^{1,1}V=\wedge^2V$.

 Let $X=G/P$ be a rational homogeneous variety, where
  $P$ is a parabolic subgroup (\cite{Ko}, \cite{F-H}). 
We fix a splitting
${Lie P}={Lie R}\oplus{Lie N}={\cal R}\oplus { \cal N}$, where
${\cal R}$ is reductive and ${\cal N}$ is the nilpotent radical. 
A representation of $P$ is completely reducible iff it is trivial
 on $N$ (see \cite{Ise} or \cite{Ot}).
In this case the representations are determined by their restriction
on $R$.

\smallskip

{\bf Homogeneus vector bundles} The group $G$ is a principal bundle over 
$X=G/P$
with fiber $P$. Denote by $z$ the point of $X$ which is fixed by
$P$, corresponding to the lateral class $P\in G/P$. Any $G$-homogeneous
vector bundle $E$ with fiber $E(z)$ over $z$ is induced by this principal
bundle via
a representation $\rho\colon P\to GL(E(z))$.
We denote $E=E_{[\rho]}$. 
 Equivalently,
$E_{[\rho]}$ can be defined as the
quotient $G\times_{\rho}E(z)$ of $G\times E(z)$ via
the equivalence relation $\sim$ where
$(g,v)\sim (g',v')$ iff there exists $p\in P$ such that $g=g'p$
and $v=\rho(p^{-1})v'$.

 We denote by $E_{\lambda}$
the homogeneous bundle corresponding to the irreducible 
representation of $P$ with maximal weight $\lambda$.
Here $\lambda$ belongs to the fundamental Weyl chamber of the reductive part 
of $P$ (see the beginning os \S 4).

\smallskip

{\bf Hermitian symmetric varieties}
 We recall that the tangent bundle of $X$ is defined by the
 adjoint representation
 over $Lie G/Lie P$.
According to Kostant, we say $X$ is a Hermitian symmetric variety if
 the above adjoint representation is trivial
on $N$. This is equivalent to ask $[{ \cal N},{ \cal N}]=0$.
The Hermitian symmetric varieties were classified by Cartan and their list 
is well known. They are product of irreducible ones. 
 The irreducible ones are grassmannians, quadrics, spinor varieties, 
maximal lagrangian grassmannians and two  varieties of exceptional type of dimension 
$16$ and $27$ (see \thmref{table} for the
precise list). For a modern treatment see
\cite{Ko} or \cite{L-M}. 
According to the corresponding Dynkin diagram, an irreducible Hermitian symmetric variety
is called of type ADE if $G=SL(m), Spin(2m), E_6$ or $E_7$. Only odd quadrics and maximal lagrangian grassmannians
are left, which are called of type BC. A Hermitian symmetric variety is called of type ADE
if it is the product of irreducible Hermitian symmetric varieties of type ADE.
Hermitian symmetric varieties of ADE type have two interesting features that we want to underline. The first one is that
when $X\subset P(V)$ is a minimal homogeneous embedding, then $V$ has a weight structure which make it isomorphic 
to the cohomology group $H^*(X, {\bf C})$. The second one is that the degree of Schubert cycles can be computed easily from the
 Hasse quiver, as in next paragraph. The reason why we have to restrict to the ADE type in the computation
of cohomology is explained in \lemref{needADE} and \lemref{needADE2}.
In all the irreducible cases we have $Pic(X)={\bf Z}$. Thus on irreducible Hermitian symmetric varieties the first 
Chern class $c_1(E)$ of a bundle $E$ can be identified with an integer, 
and the slope is by definition $\mu(E)=c_1(E)/rk(E)\in{\bf Q}$.
On any Hermitian symmetric variety $X=X_1\times\ldots\times X_r$ where $X_i$ are irreducible
there are several possible choices of slopes. With obvious notations, if $c_1(E)=(c_1^1,\ldots, c_1^r)\in{\bf Z}^r$
and $a=(a_1,\ldots, a_r)\in {\bf Q}^r$ then we define
$\mu_a(E)=\frac{\sum c_1^ia_i}{rk(E)}\in{\bf Q}$. 

It is easy to check (see e.g. \cite{Ram} 5.2)
that
$\mu_a(E_{\sum n_i\lambda_i})=\sum n_i\mu_a(E_{\lambda_i})$.

\smallskip

{\bf The Hasse quiver}
Quivers will be recalled in \S 5. For this paragraph it is enough to know that a quiver is just
an oriented graph.
If $X$ is a rational homogeneous variety, the cohomology $H^*(X, {\bf Z})$ can be organized
in a quiver in the following way. Consider the action of a Borel subgroup $B\subset P$ on $X$. Then it is well known that
$X$ is divided in a finite union of orbits, their closures are called the Schubert celles and form an additive basis
 $H^*(X, {\bf Z})$. The vertices of the Hasse quiver ${\cal H}_X$ are the Schubert celles, we draw an arrow
 between $X_{\omega}\in H^{2p}(X, {\bf Z})$ and $X_{\omega'}\in H^{2p+2}(X, {\bf Z})$
 if $X_{\omega}\supset X_{\omega'}$. If $X$ is a Hermitian symmetric variety
 the additive basis of $H^{2p}(X, {\bf Z})$  corresponds to the direct summands of $\Omega^p$.
 If $X$ is Hermitian symmetric, the degrees of the Schubert cycles in the 
 homogeneous minimal embedding are computed as the number of paths in the Hasse quiver which starts
 from the corresponding vertex. We learned this fact from L. Manivel (see \cite{I-M}).

\smallskip

{\bf The filtration of a homogeneous bundle
and the functor $gr$} 
Let $E$ be a homogeneous bundle on an irreducible Hermitian symmetric variety.

We define $gr E=\oplus_i E_i/E_{i-1}$ for any
filtration $0\subset E_1\subset\ldots\subset E_k=E$ such that 
$E_i/E_{i-1}$ is completely reducible.  The  graded bundle $gr E$ does not depend on the filtration,
 in fact it is given by the restriction
of the representation giving $E$ to the reductive part $R$ of $P$.

 For example the Euler sequence on ${\bf P}={\bf P}(V)$ tells us that
$gr ({\cal O}(1)_{{\bf P}}\otimes V)= {\cal O}_{{\bf P}}\oplus T{\bf P}$.

The functor $E\mapsto gr E$ from $P$-$mod$ to $R$-$mod$
(which in the literature is often denoted as $Ind^P_R$)
is exact. It is easy to check the formulas
\[(gr E)^*=gr(E^*)\quad gr(E\oplus F)= gr E\oplus gr F\quad
gr(E\otimes F)= gr E\otimes gr F\]

{\bf The spectral sequence abutting
to the cohomology} The Borel-Weil-Bott
 theorem describes
 the cohomology of the irreducible homogeneous bundles $E$. 
It says that $H^*(E)$
is an irreducible $G$-module.
For any homogeneous bundle and for any filtration
there is a spectral sequence 
 abutting to the  cohomology of the bundle.
Precisely, if $gr E=\oplus_{i=1}^kA_i$ as before, we have 
$E^1_{p,q}=H^{p+q}(A_{k-p})$ abutting
to $E^{\infty}_{p,q}$ where $H^i(E)=\oplus_{p+q=i}E^{\infty}_{p,q}$. 
\thmref{filt1} will give a more efficient way to compute
$H^i(E)$.

\smallskip

{\bf Yoneda product} We recall the  Yoneda product
 on $Ext$ according
 to \cite{Ei}
exerc. A3.27.
For any homogeneous bundles $E$, $F$ and $K$ there is an equivariant
Yoneda product
\[Ext^i(E,F)\otimes Ext^j(F,K)\to Ext^{i+j}(E,K)\]
and this product is associative.
In particular in the case $E=F=K$ and $i=j=1$
we get a (non symmetric) bilinear map, whose symmetric part induces
 a quadratic morphism
\[Ext^1(E,E)\to Ext^2(E,E)\]
In particular, since it preserves the invariant part, it gives
\[m\colon Ext^1(E,E)^G\to Ext^2(E,E)^G\]

{\bf Tensor product of two irreducible representations}
 Let $\lambda$ and $\nu$ two weights in the fundamental Weyl
 chamber of a Lie algebra $K$. The tensor product of the corresponding 
representations $V_{\lambda}\otimes V_{\nu}$ can be expressed as a 
sum $\oplus c_{\lambda\nu\kappa}V_{\kappa}$ where $c_{\lambda\nu\kappa}$
are integers (counting the multiplicities).
When $K=Lie SL(n)$ the integers $c_{\lambda\nu\kappa}$
can be computed by the so called Littlewood-Richardson rule
(see \cite{F-H}). 
A more conceptual algorithm was later conjectured by Weyman and proved by 
Littelmann in \cite{Li}; this algorithm holds
 for an arbitrary simple Lie groups. 
Let ${\nu_1}=\nu,{\nu}_2,\ldots ,{\nu}_k  $
be all the weights of  $V_{\nu}$.   
Littelmann proves that
\begin{equation}
\label{littelmann}
V_{\lambda}\otimes V_{\nu}=\oplus_{i\in I}V_{\lambda+\nu_i}
\end{equation}
where $I$ is a subset of $\{1,\ldots ,k\}$ such that
 the weights $\nu_i$ for $i\in I$ correspond exactly to the standard Young 
tableaux of the form corresponding to $\nu$ which are $\lambda$-dominant
(see \cite{Li} for the precise definitions).
A particular interesting case is when
$\lambda+\nu_i$ are all dominant for
$i=1,\ldots ,k$, this is true when $\lambda\gg 0$. In this case we have the whole decomposition
\[V_{\lambda}\otimes V_{\nu}=\oplus_{i=1}^kV_{\lambda+\nu_i}\]
(see also \cite{F-H} exerc. 25.33).
Formula (\ref{littelmann}) above applied to vector bundles gives
$$E_{\lambda}\otimes E_{\nu}=\oplus_{i\in I}E_{\lambda+\nu_i} $$
where all the direct summand in the right side have the same slope
(see \cite{Ram} or \cite{Ot}).

\section {$P$-$mod$ and the category of Higgs bundles}

Let $X$ be a Hermitian symmetric variety.
We recall that ${\cal N}$ is an ${\cal R}$-module
with the adjoint action. Our starting point is the following
\begin{thm0}
\label{starting}

(i) Given a ${\cal P}$-module $E$ on $X$, the action of ${\cal N}$
over $E$ induces a morphism of ${\cal R}$-modules
\[\theta\colon{\cal N}\otimes gr E\rig{}gr E\]
such that $\theta\wedge\theta=0$
in $Hom(\wedge^2{\cal N}\otimes gr E,gr E)$

(ii) Conversely given an ${\cal R}$-module $F$ on $X$
and a morphism of ${\cal R}$-modules
\[\theta\colon{\cal N}\otimes F\rig{}F\]
such that $\theta\wedge\theta=0$
then $\theta$ extends uniquely to an action of ${\cal P}$ over $F$,
giving a bundle $E$ such that $gr E=F$.
\end{thm0}

{\it Proof} (i) For every $r\in {\cal R}$,  $n\in {\cal N}$,
 $f\in F$, since $E$ is a ${\cal P}$-module
we have
\[r\cdot(n\cdot f)=n\cdot(r\cdot f)+[r,n]\cdot f\]
that is 
\[r\cdot(\theta(n\otimes f))=
\theta (n\otimes (r\cdot f))+\theta([r,n]\otimes f)=
\theta(r\cdot(n\otimes f)) \]
so that $\theta$ is ${\cal R}$-equivariant.
Moreover for any $n_1, n_2\in {\cal N}$
\[\theta\wedge\theta\left((n_1\wedge n_2)\otimes f\right)=
n_1\cdot(n_2\cdot f)-n_2\cdot(n_1\cdot f)=[n_1,n_2]\cdot f=0\]
because $[{\cal N},{\cal N}]=0$ and this is equivalent to
$\theta\wedge\theta=0$.

(ii)  We have for any  $r+n\in{\cal R}\oplus{\cal N}={\cal P}$
$$(r+n)\cdot f:=r\cdot f+\theta(n\otimes f)$$ and we have to prove
that for any $p_1, p_2\in {\cal P}={\cal R}\oplus{\cal N}$
we have 
\begin{equation}
\label{pmodule}
[p_1,p_2]\cdot f=p_1\cdot(p_2\cdot f)-p_2\cdot(p_1\cdot f)
\end{equation}
We distinguish three cases.

If $p_1,p_2\in {\cal R}$  then (\ref{pmodule}) is true
because $F$ is an ${\cal R}$-module.

If $p_1,p_2\in {\cal N}$ then  $[p_1,p_2]=0$ and
(\ref{pmodule}) is true
because $\theta\wedge\theta=0$.

If $p_1\in {\cal R}$, $p_2\in {\cal N}$ we have $[p_1,p_2]\in{\cal N}$
and
\[[p_1,p_2]\cdot f+p_2\cdot(p_1\cdot f)=
\theta(p_1\cdot(p_2\otimes f))=p_1\theta(p_2\otimes f)=
p_1\cdot(p_2\cdot f)\]
because  $\theta$ is ${\cal R}$-equivariant.
\qedd

\medskip

\thmref{starting} allows one
to construct a ${\cal P}$-module in two steps:
the first step is to give the ${\cal R}$-module $F$,
which encodes the discrete part, the second step
is to give $\theta$, which encodes the continuous part.
This will be made precise in \S 7 about moduli spaces. At present
it is convenient to reformulate 
\thmref{starting} in terms of vector bundles.

We have seen in the introduction
that on a Hermitian symmetric variety the
${\cal P}$-module ${\cal N}$ corresponds to $\Omega^1_X$.
Since $[{\cal N},{\cal N}]=0$, $\Omega^1_X$ is completely reducible.
Let $E$ be a $G$-homogeneous bundle $E$. 
The action of $\cal N$ over the $\cal R$-module
$gr E$ induces by
\thmref{starting} an ${\cal R}$-equivariant morphism
of completely reducible representations
${\cal N}\otimes gr E\to gr E$, hence we get a $G$-equivariant morphism
$\theta\in Hom(gr E, gr E\otimes T_X)^G$ such that
$\theta\wedge\theta=0$ .
To any $E$ we can associate the pair
$(gr E,\theta)$. Such pairs are analogous
to what is called in \cite{Simp} a Higgs bundle.
The pairs $(gr E,\theta)$ are
the natural extension of Higgs bundles for rational homogeneous 
varieties, where $T_X$ is globally generated, 
so we maintain the terminology of Higgs bundles also in this case.
             
More precisely, we have
\begin{defn0} Let $X$ be a Hermitian symmetric variety. A Higgs bundle
on $X$ is a pair $(F,\theta)$ where $F$ is an ${R}$-module
and $\theta\colon F\rig{}F\otimes T_X$ 
is $G$-equivariant and satisfies
$\theta\wedge\theta=0$. 
\end{defn0}

Higgs bundles form an abelian category, where a morphism between two Higgs bundles
$(F_1,\theta_1)$ and $(F_2,\theta_2)$ is a $G$-equivariant
morphism $f\colon F_1\rig{}F_2$ such that
$(f\otimes id)\theta_1=\theta_2f$.
Hence \thmref{starting}
can be reformulated in the following way:

\begin{thm0}
\label{starting2}   
Let $X=G/P$ be a Hermitian symmetric variety.
 There is an equivalence of categories between

(i) $G$-homogeneous bundles over $X$ 

(ii) Higgs bundles $(F,\theta)$ over $X$ 
\end{thm0}

{\bf Remark} On any rational homogeneous variety,
the category of $G$-homogeneous bundle is equivalent
to the category of pairs $(F,\theta)$ where $F$ is an ${R}$-module
and $\theta\colon F  \rig{}F\otimes T_X$ 
is $G$-equivariant and satisfies certain relations.

\section{ The Borel-Weil-Bott theorem}

It is well known that the hyperplanes orthogonal to the roots of $G$ divide ${\cal H}^{\vee}$ into regions called Weyl chambers. 
The fundamental Weyl chamber $D$ of $G$ is
$$D=\{\sum x_i\lambda_i |x_i\ge 0\}$$
and it contains exactly the dominant weights.
 The Weyl group ${ W}$ acts 
in simple transitive way 
as a group of isometries
on the Weyl chambers.  Following \cite{Ko} we denote
$g=\sum\lambda_i$. 
Any homogeneous variety with $Pic={\bf Z}$ is  the quotient
$X=G/P(\alpha_j)$ for some $j$,
where the Lie algebra of $P(\alpha_j)$ is spanned
by the Cartan subalgebra, by the eigenspaces of the negative roots
 and by the eigenspaces of the positive roots 
$\alpha=\sum n_i\alpha_i$
such that $n_i\ge 0$ for any $i$ and $n_j=0$.

The reductive part of $P(\alpha_j)$ has its own fundamental Weyl chamber 
$D_1\supset~D$ defined by 
$$D_1=\{\sum x_i\lambda_i |x_i\ge 0 \hbox{\ for } i\neq j\}$$ 
$D_1$ contains exactly the maximal weights of the irreducible representations
of $P(\alpha_j)$.
Let \[W^1=\{ w\in W| wD\subset D_1 \}\] 
(see \cite{Ko} Remark 5.13).
The cardinality of $W^1$ divides the order of $W$. 

Let $H_{\phi}$ be the hyperplane orthogonal to the root $\phi$ and
 $r_{\phi}$ be the reflection with respect to $H_{\phi}$.
It is well known that the reflections $r_{\alpha_i}$ generate the Weyl group. 

Let $Y_{\phi}=H_{\phi}-g$.

Let $\xi_1,\ldots , \xi_m$ be the weights
of the representation giving the bundle
$\Omega^1_X$, where $m=\dim X$.
Let $s_j$ for $j=1,\ldots ,m$ be the reflection through
$Y_{\xi_j}$; note that for any weight $\lambda$
\begin{equation}
s_j(\lambda)=r_{\xi_j}(\lambda +g)-g
\end{equation}
 thus
$s_j$ and $r_{\xi_j}$ are conjugate elements
in $Iso({\cal H}^{\vee})$.
It follows that if $w=r_{\xi_1}\cdot\ldots\cdot r_{\xi_p}$
then $w(\lambda+g)-g=s_1\cdot\ldots\cdot s_p(\lambda)$.

An element $\nu\in Z$ is called regular if $(\nu,\phi)\neq 0$
for any root $\phi$, otherwise it is called singular. Observe that
$\nu$ is singular iff $\nu\in H_{\phi}$ for some root $\phi$.

Denote (see \cite{Ko} Remark 6.4)
\[D_1^0=\{\xi\in D_1 | g+\xi\hbox{\ is regular}\}\]

$D_1^0$ consists in the subset of $D_1$ obtained removing
exactly the $Y_{\xi_j}$. Hence a convenient composition of $s_j$
brings $D$ into the several "chambers" in which
$D_1^0$ is divided, which we call {\it Bott chambers}
 (do not confuse them with the
usual Weyl chambers). The Bott chambers  are obtained
by performing a slight "separation" on the Weyl chambers,
see the following picture in the case of ${\bf P}^2=SL(3)/P(\alpha_1)$
where the three Bott chambers are shadowed

\bigskip

\hspace*{2cm}
\includegraphics[scale=0.3]{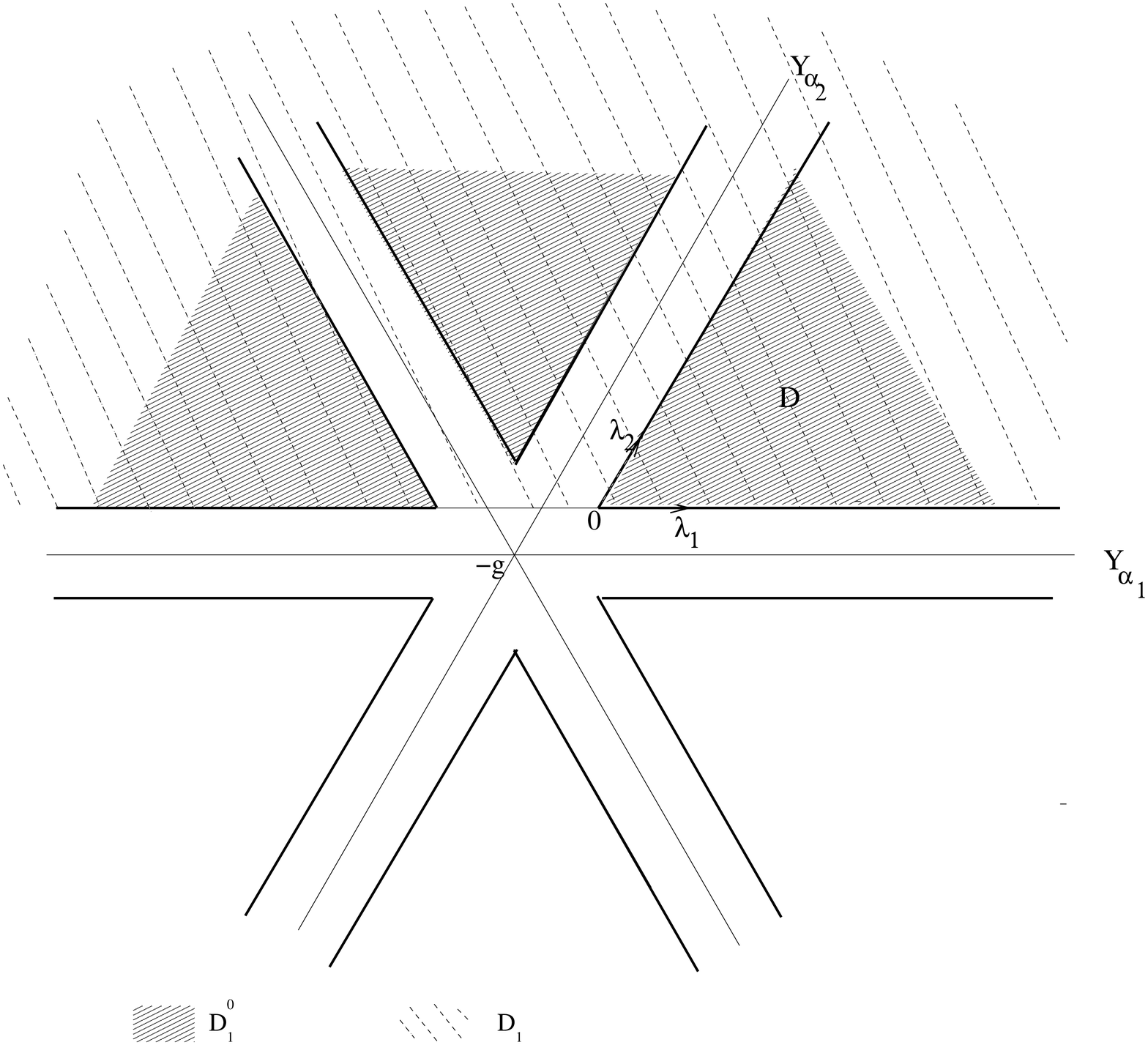}
\bigskip

Now for any  $w\in W$ the {\it length} $l(w)$ is defined as the 
minimum number of reflections $r_{\alpha}$
(with $\alpha$ root) needed to obtain $w$.
Any Bott chamber has its own length. Two Bott chambers are said to be adjacent if
they have a common hyperplane in their boundary. The lengths of two
Bott chambers are consecutive integers.

We state the Bott Theorem 
(compare with \cite{Ko} theorem 5.14)

\begin{thm0}
\label{bottclassico} {\bf (Bott)}
If $\lambda\in D_1$, then $\exists !$
$w\in W$ s.t. $w^{-1}\in W^1$ and
$ w(\lambda+g) \in D$. 

(i) If
$w(\lambda+g) $ belongs to the interior of $D$ then setting
$\nu=w(\lambda+g)-g$ we have 
$H^{l(w)}(E_\lambda)=V_{\nu}$ and $H^j(E_\lambda)=0$ 
for $j\neq l(w)$.

In particular if $\lambda\in D$ (thus $w$ is the identity) then
$H^0(E_\lambda)=V_\lambda$ and $H^i(E_\lambda)=0$ for $i>0$.

(ii) If
$w(\lambda+g) $ belongs to the boundary of $D$ then 
$H^j(E_\lambda)=0\quad\forall j$. 
\end{thm0}

We recall the  result of Kostant
 (\cite{Ko} Corollary 8.2):
\[ \#\{w\in W^1 | l(w)=i\} = \dim H^{2i}(X, {\bf C})\]

in particular 
\begin{equation}\label{chi}\# W^1 =\chi (X,  {\bf C})
\end{equation}

We explain now the relation of the previous result with the Bott
theorem. By Hodge-Deligne theory 
$H^{2i}(X, {\bf C})$ is isomorphic to
$H^i(X, \Omega_X^i)=H^i(X, \Omega_X^i)^G$. Moreover
for any irreducible Hermitian symmetric varieties the bundle 
$\Omega^1$
is irreducible and $\Omega^i$ splits as a sum
of direct summands and the number of these summands is equal to
$\dim H^{2i}(X, {\bf C})$. Moreover on $X=X_1\times\ldots\times X_r$
with projections $p_i$ we have $\Omega^1_X=\oplus p_i^*\Omega^1_{X_i}$.
The vertices $\lambda$ of the Bott chambers correspond
exactly to the direct summands of $\Omega^i$ for some $i$.
Indeed for any such a vertex $\lambda$ there exists
$w$ as in the Bott theorem (thus $w^{-1}\in W^1$)
such that $w(\lambda+g)-g=0$ (i.e. $\lambda=w^{-1}(g)-g$)
and  $l(w)=i$.

We note the following consequences of the results of Bott and Kostant

\begin{coro0}
\label{hjhom}
Let $E$ be a completely reducible bundle on $X$ 
Hermitian symmetric variety.
Then $H^j(E)^G$ is isomorphic to
$Hom(\Omega^j,E)^G$. This means that, when $E$
is irreducible, $H^j(E)^G\neq 0$
if and only if $E$ is a direct summand of $\Omega^j$.
\end{coro0}

{\it Proof} We may suppose $E=E_{\lambda}$. We have
$Hom(\Omega^j,E_{\lambda})^G\neq 0$ iff
$E_{\lambda}$ is a direct summand of $\Omega^j$ and
in this case it is isomorphic to ${\bf C}$. By Bott theorem we have 
$H^j (E_{\lambda})^G\neq 0$ iff $H^j (E_{\lambda})={\bf C}$
and this is true only if 
$w(\lambda+g)-g=0$ ($w$ as in the Bott theorem) and $l(w)=j$. 
These cases are exactly when $E_{\lambda}$ is a direct summand
of $\Omega^j$.\qedd

\begin{thm0}
\label{homext}
Let $X=G/P$ be a Hermitian symmetric variety.

(i) There is a natural isomorphism
$Hom( E_{\lambda}\otimes\Omega_X^i, E_{\nu})^G
\to Ext^i(E_{\lambda}, E_{\nu})^G$ $\forall\lambda, \nu\in D_1$.
Both spaces are isomorphic to ${\bf C}$ or to $0$
for $i=1$.

(ii) If $X$ is irreducible and $Ext^i(E_{\lambda}, E_{\nu})^G\neq 0$ then 
$\mu(E_{\nu})=\mu(E_{\lambda})+i\mu(\Omega^1)$.

(iii) If $X=X_1\times\ldots\times X_r$, product of irreducible ones, and $Ext^i(E_{\lambda}, E_{\nu})^G\neq 0$
define $a_i=\frac{1}{\mu(\Omega^1_{X_i})}$.
Then  with this choice for any $i$ we have
$\mu_a(\Omega^1_X)=\mu_a(p_i^*\Omega^1_{X_i})=1$ (see \S 2) and we get
$\mu_a(E_{\nu})=\mu_a(E_{\lambda})+i$.

\end{thm0}

{\it Proof:} (i) By
\corref{hjhom} only the last statement  needs an explanation.
In fact
all the irreducible components of
$E_{\lambda}\otimes\Omega_X^1$ have multiplicity one.
Indeed look at (\ref{littelmann})
and observe that eigenspaces of the roots of $G$ 
have dimension $1$. 

(ii) All direct summands of 
$E_{\lambda}\otimes\Omega_X^i$
have the same $\mu$ equal to $\mu(E_{\lambda} \otimes \Omega_X^i)=
\mu(E_{\lambda})+i\mu(\Omega_X^1)$.

(iii) follows immediately as in (ii).
\qedd
\medskip

{\bf Remark} For $i\ge 2$
there are some irreducible components
of $E_{\lambda}\otimes\Omega_X^i$ which appear with multiplicity
$\ge 2$. For example in the Grassmannian
$Gr({\bf P}^1, {\bf P}^3)=SL(4)/P(\alpha_2)$, 
let $T$ be the tangent bundle. We have
that $Ext^2(T, T(-2))^G$ contains $H^2(\Omega^2)={\bf C}^2$,
and correspondingly $T\otimes \Omega^2$ contains two copies of
$T(-2)$ . Indeed $\Omega^2$ splits
into two irreducible summands and there is a copy of
$T(-2)$ for each of these summands. In the case of quadrics
$Q_n$ with $n\ge 5$, the list of weights of the irreducible 
$\Omega^2$ contains a weight of multiplicity 
$[{\frac{n-1}{2}}]$, in this case
for $\lambda\gg 0$ the tensor product $E_{\lambda}\otimes\Omega^2$
contains a direct summand with multiplicity $[{\frac{n-1}{2}}]$.
 In the case $X={\bf P}^n$ all irreducible summands of
 $E_{\lambda}\otimes\Omega^2$ appear with multiplicity one by
 the formula (\cite{F-H} (6.9)), indeed in this case all the
 weights of $\Omega^2$ are distinct. 
 
\medskip

\begin{coro0}
\label{irrexti}
If $E$ is an irreducible bundle on  a Hermitian symmetric 
variety then $Ext^i(E,E)^G=0$ for $i>0$.
\end{coro0}

{\it Proof} Apply \thmref{homext} for $\lambda=\nu$.
\qedd

\begin{coro0}  
\label{lambdalambda'}
For every $i<\dim X$ and $\lambda\in D_1$ there are
$\lambda'$ and $s_j$ such that $\lambda'=s_j(\lambda)$
and $H^i(E_{\lambda})=H^{i+1}(E_{\lambda'})$ or
$H^i(E_{\lambda})=H^{i-1}(E_{\lambda'})$.
In particular $\lambda$ and $\lambda'$ differ by a multiple of
$\xi_j$. There is exactly
one of such $\lambda'$ in every Bott chamber
having a common boundary with the chamber containing
$E_{\lambda}$.  
\end{coro0}

{\it Proof} Consider the vertex $\lambda_0$
of the Bott chamber containing $\lambda$.
Then consider all the $s_j$ such that  $s_j(\lambda_0)$
is the maximal weight of a summand of $\Omega^{i+1}$.
Such $s_j$'s work. 

\qedd   

\medskip

{\bf Remark}  $\lambda'$ and $s_j$  of the previous corollary
are unique in the case of ${\bf P}^n$, but they are not unique for 
general Grassmannians.

In the following tables we list all the vertices of the Bott
chambers in the cases ${\bf P}^4$ and $Gr(1,4)$. The $4$-ple
$(x_1,x_2,x_3,x_4)$ denotes
the weight $\sum x_i\lambda_i$. An arrow labelled with the root
$\beta$ means the reflection
$$\cdot \mapsto r_{\beta}(\cdot +g)-g$$
So the arrow labelled with $-\xi_j$ means the reflection $s_j$.
For example $(-2,1,0,0)=r_{\alpha_1}\left((0,0,0,0)+g\right)-g=s_1(0,0,0,0)$.
To check the tables it can be useful \lemref{alphapn}.

\[{\bf P}^4\]
\[\begin{array}{l}
\bullet\ (0,0,0,0)\cr
\dow{\alpha_1}\cr
\bullet\ (-2,1,0,0)\cr
 \dow{\alpha_1+\alpha_2}\cr
\bullet\ (-3,0,1,0)\cr
  \dow{\alpha_1+\alpha_2+\alpha_3}\cr
\bullet\ (-4,0,0,1)\cr
   \dow{\alpha_1+\alpha_2+\alpha_3+\alpha_4}\cr
\bullet\ (-5,0,0,0)\cr
\end{array}\]
 \vfill\eject
\[{\bf Gr(1,4)}\]
\[\begin{array}{rllrl}
&&\bullet\rlap{(0,0,0,0)}\cr
&&\dow{\alpha_2}\cr
&&\bullet\rlap{(1,-2,1,0)}\cr
&\lswar{\alpha_1+\alpha_2\qquad }&&\sear{\alpha_2+\alpha_3}\cr
(0,-3,2,0)\bullet &&&&\bullet\ (2,-3,0,1)\cr
\ldow{\alpha_2+\alpha_3\ }&&\swar{\alpha_1+\alpha_2}&&\dow{\alpha_2+\alpha_3+\alpha_4}\cr
(1,-4,1,1)\bullet &&&&\bullet\ (3,-4,0,0)\cr
\ldow{\alpha_1+\alpha_2+\alpha_3\ }&&\sear{\alpha_2+\alpha_3+\alpha_4}
&&\dow{\ \alpha_1+\alpha_2}\cr
(0,-4,0,2)\bullet &&&&\bullet\ (2,-5,1,0)\cr
&\lsear{\alpha_2+\alpha_3+\alpha_4\qquad }&&\swar{\alpha_1+\alpha_2+\alpha_3}\cr
&&\bullet\rlap{(1,-5,0,1)}\cr
&&\dow{\alpha_1+\alpha_2+\alpha_3+\alpha_4}\cr
&&\bullet\rlap{(0,-5,0,0)}\cr
\end{array}\]

Of course the above graphs are exactly the Hasse quivers
${\cal H}_{{\bf P}^4}$ and ${\cal H}_{Gr(1,4)}$.

On ${\bf P}^n$ we have a simplification of the Bott theorem.
In this case $\Omega^p$
are irreducible $\forall p$.

\begin{prop0} 
\label{bottpn}{\bf (Bott on $ {\bf P}^n $)}
Let $X={\bf P}^n =SL(n+1)/P(\alpha_1)$
 
(i) if $\lambda$ is any weight and $\exists i\in{\bf N}$ s.t. 
$\nu:=r_{\alpha_i}\ldots r_{\alpha_{1}}(\lambda+g)-g\in D$ then
$H^i(E_\lambda)=V_{\nu}$ and $H^j(E_\lambda)=0$ for $j\neq i$.

In particular, if $\lambda\in D$ then
$H^0(E_\lambda)=V_\lambda$ and $H^i(E_\lambda)=0$ for $i>0$.

(ii) in the remaining cases $H^j(E_\lambda)=0\quad\forall j$ 
\end{prop0}

{\it Proof } 
 It is sufficient, by \thmref{bottclassico},
 to prove that
$W^1=\{ r_{\alpha_1}\cdot\ldots\cdot r_{\alpha_{i}} |
i\in\{1\ldots n\}\}\cup \{1\}$. 
It is well known that
$r_{\alpha_i}(\lambda_j)$ is equal to
$\lambda_j$ if $j\neq i$, and to
$\lambda_{j-1}-\lambda_j+\lambda_{j+1}$
if $j=i$ (with the convention that $\lambda_0=\lambda_{n+1}=0$).
It holds that
 $$r_{\alpha_1}\cdot\ldots\cdot r_{\alpha_{i}}
(\sum_{j=1}^np_j\lambda_j)=(-\sum_{j=1}^ip_i)\lambda_1+
\sum_{j=1}^ip_j\lambda_{j+1}+\sum_{j=i+1}^np_j\lambda_{j}$$
(to check it prove that 
$r_{\alpha_1}\ldots r_{\alpha_{i}} (\lambda_j)$ is equal to
$\lambda_j$ if $j>i$ and it is equal to 
$r_{\alpha_1}\ldots r_{\alpha_{j}} (\lambda_j)=
-\lambda_1+\lambda_{j+1}$ if $j\le i$).

Hence the elements $r_{\alpha_1}\cdot\ldots\cdot r_{\alpha_{i}}$
belong to $W^1$ for $i=1$ to $n$, so these elements, together with the identity, fill $W^1$ by
 (\ref{chi}). The last remark is that
$(r_{\alpha_1}\ldots r_{\alpha_{i}})^{-1}=
r_{\alpha_i}\ldots r_{\alpha_{1}} $
\qedd
\medskip

The point (iv) of the following lemma gives an alternative way
to express point (i) of the Bott theorem.

\begin{lemma0}\label{alphapn}
On $ {\bf P}^n $ we have for $i=1\ldots n$
\[\begin{array}{l}
(i)\ \xi_i=-\alpha_1+\ldots-\alpha_i\cr
(ii)\ \alpha_{1}+\ldots +\alpha_{i+1}=
\left(r_{\alpha_1}\ldots r_{ \alpha_{i}}\right)(\alpha_{i+1})\cr
(iii)\ r_{\xi_{i+1}}=
(r_{\alpha_{i}}\cdot\ldots\cdot r_{\alpha_1})^{-1}
r_{\alpha_{i+1}}(r_{\alpha_i}\cdot\ldots\cdot r_{\alpha_1})\cr
(iv)\ r_{\xi_1}\cdot\ldots\cdot r_{\xi_i}=
r_{\alpha_i}\cdot\ldots\cdot r_{\alpha_1}\cr
\end{array}\]
\end{lemma0}

{\it Proof} Straightforward (for (iii) observe that
by (ii) $r_{\alpha_1+\ldots +\alpha_{i+1}}
=r_{r_{\alpha_1}\ldots r_{ \alpha_{i}}(\alpha_{i+1})}$).
 
\qedd

\begin{coro0}  On ${\bf P}^n$
if 
$\lambda=s_{i+1}(\lambda')$ then
$H^i(E_{\lambda})=H^{i+1}(E_{\lambda'})$  . The converse holds
if $H^i(E_{\lambda})\neq 0$.

In particular $\lambda$ and $\lambda'$ differ by a multiple of
$\alpha_1+\ldots +\alpha_{i+1}$. Precisely
if $\lambda= \sum_{j=1}^np_j\lambda_j$ then
$\lambda'-\lambda=
-\sum_{j=1}^{i+1}(p_j+1)(\alpha_1+\ldots +\alpha_{i+1})$.

\end{coro0}

{\it Proof} By \propref{bottpn} only the converse needs
to be proved.  
If $H^{i+1}(E_{\lambda'})=H^i(E_{\lambda})\neq 0$ then 
by the Bott theorem $h(\lambda+g)=r_{ \alpha_{i+1} } h(\lambda'+g)$
where $h=r_{\alpha_i}\ldots r_{\alpha_{1}} $
and this implies that
$h(\lambda+g)-h(\lambda'+g)$ is parallel to
$\alpha_{i+1} $, that is $\lambda-\lambda'$
is parallel to 
$h^{-1}\alpha_{i+1}=\alpha_{1}+\ldots +\alpha_{i+1} $
(by \lemref{alphapn} (ii)). Moreover the last formula holds because
 $(\lambda+g,\xi_{i+1})=-\sum_{j=1}^{i+1}(p_{j}+1)$.
\qedd

\section{The Quiver and its Relations}

For a quick introduction to theory of quivers and their representations 
we refer to \cite{King}. More details about quivers with relations can be found
in \cite{G-R} or in \cite{Hi1}.
 
\begin{defn0}
\label{quiverdef}
A {\bf quiver} is an oriented  graph ${\cal Q}$ with the set ${\cal Q}_0$ of 
 points and the set  ${\cal Q}_1$ of arrows. There are two maps
$h,t\colon{\cal Q}_1\to {\cal Q}_0$ which indicate respectively
the head (sink) and the tail (source) of each arrow.

A {\bf path} in ${\cal Q}$ is  a formal composition of arrows $ \beta_m ...
\beta_1$ where the tail of an arrow is the head of the previous one.
Paths can be summed and composed in natural way, defining the {\bf path algebra} 
${\bf C}{\cal Q}$. It is graded by pairs in ${\cal Q}_0$.

A {\bf relation} in ${\cal Q}$ is a linear form $\lambda_1 c_1+...+
\lambda_m c_m$ where $c_i$ are paths in ${\cal Q}$ with a common tail and
 a common head and $\lambda_i \in {\bf C}$.

A {\bf representation of a quiver} ${\cal Q} =({\cal Q}_0, 
{\cal Q}_1)$ is the couple of a set of vector 
spaces $\{X_i\}_{i \in {\cal Q}_0} $ and of a set of linear maps 
$ \{\varphi_{\beta} \}_{\beta \in {\cal Q}_1}$ where  $\varphi_{\beta} : X_i 
\rightarrow X_j$ if $\beta$ is an arrow from $i$ to $j$.

Let ${\cal R}$ be a homogeneous ideal in the path algebra.
A {\bf representation of a  quiver ${\cal Q}$ with relations ${\cal R}$} 
is a representation of the  quiver s.t. 
$$\sum_j \lambda_j \varphi^j_{1} ...\varphi^j_{m_j}=0$$ 
for every $ \sum_j \lambda_j \beta^j_{1} ...\beta^j_{m_j} \in {\cal R}$.

Let $(X_i, \varphi_{\beta})_{i \in {\cal Q}_0,\; \beta \in {\cal Q}_1}$ and 
$(Y_i, \psi_{\beta})_{i \in {\cal Q}_0,\; \beta \in {\cal Q}_1}$ be two 
representations of the quiver ${\cal Q}= ({\cal Q}_0 , {\cal Q}_1)$.
A {\bf morphism} $f $ from $(X_i, \varphi_{\beta})_{i \in {\cal Q}_0,\; 
\beta \in {\cal Q}_1}$ to
$(Y_i, \psi_{\beta})_{i \in {\cal Q}_0,\; \beta \in {\cal Q}_1}$ is a
 set of linear maps $f_i : X_i \rightarrow Y_i $, $i \in {\cal Q}_0$ s.t. for
 every  $\beta \in {\cal Q}_1$,  $\beta$  arrow from $ i$ to 
 $j$, the following diagram is commutative: 
$$\begin{array}{ccc} X_i   
 & \stackrel{f_i}{\longrightarrow}
 &  Y_i \\ 
\varphi_{\beta} \downarrow & &  \downarrow \psi_{\beta}
\\   X_j 
 &  \stackrel{f_j}{\longrightarrow} &  Y_j
 \end{array}$$ 

It is well known (and easy to be proved) that the category of representations
of ${\cal Q}$ with relations ${\cal R}$ is equivalent to the category of
${\bf C}{\cal Q}/{\cal R}$-modules.

A quiver ${\cal Q}$ is called {\bf levelled} if there exists a function $s\colon{\cal Q}_0 \to{\bf Q}$
such that for any arrow $i\rig{}j$ we have $s(i)=s(j)+1$

\end{defn0}

Let $X=G/P$ be a Hermitian symmetric variety.
 In order to describe
all $G$-homogeneous bundles on $X$ we define a quiver ${\cal Q}_X$.

\begin{defn0}
Let ${\cal Q}_X$ be the following quiver.
The points of ${\cal Q}_X$ are the irreducible representations of 
$R$, which we
identify with irreducible $G$-homogeneous bundles over $X=G/P$, or
 with the corresponding
elements in ${\cal H}^{\vee}$.
Let $E_{\lambda}$ and $E_{\mu}$ be irreducible representations
with maximal weights  $\lambda, \mu\in D_1$.
There is an arrow in ${\cal Q}_X$ from $E_{\lambda}$ to $E_{\mu}$ iff
$Ext^1(E_{\lambda}, E_{\mu})^G\neq 0$.  The ideal of relations in 
${\cal Q}_X$ will be defined in \defref{rel12}.
\end{defn0}

Observe that if $Ext^1(E_{\lambda}, E_{\mu})^G\neq 0$ then
this group is isomorphic to ${\bf C}$, by \thmref{homext}.

\begin{coro0}
\label{levelled}
 If there is an arrow from $E_{\lambda}$ to $E_{\mu}$
then $\mu(E_{\mu})=\mu(E_{\lambda})+\mu(\Omega^1)$.
In particular the quiver is levelled (see \defref{quiverdef}) 
by $\mu_a$ of \thmref{homext} (iii)  (see \cite{Hi1}, \cite{Hi2}).

\end{coro0}

{\it Proof} By \thmref{homext}. \qedd

\medskip

\begin{coro0}
\label{negroot}
The arrows (modulo translation) between elements of the quiver can
be identified with the weights of $\Omega^1$
(which are negative roots).

\end{coro0}

{\it Proof}  From (\ref{littelmann}) it follows
$E_{\lambda}\otimes\Omega^1\subset\oplus E_{\lambda+\xi_i}$,
then we conclude by \thmref{homext}.
 \qedd

\medskip
 We postpone the description of the relations in the quiver after we have
defined the representation associated to a bundle.

\begin{defn0}
\label{defmain}
We associate to a $G$-homogeneous bundle $E$ the
following {\bf representation}
 of ${\cal Q}_X$. Let $gr E=\oplus_{\lambda}E_{\lambda}\otimes 
V_{\lambda}$, where $V_{\lambda}={\bf C}^k$ and $k$ is the number of times $E_{\lambda}$ occurs.

To the point $\lambda$ we associate the vector space $V_{\lambda}$.

For any $\lambda\in{\cal Q}_0$ let us fix a maximal vector
$v_{\lambda}\in E_{\lambda}$. For any $\xi_i$ root
of ${\cal N}$ let us fix an eigenvector $n_i\in{\cal N}$.
We have 
\begin{equation}
\label{sumtensor}Ext^1(gr E, gr E)=\oplus_{\lambda, \mu}
 Hom(V_{\lambda},V_{\mu})
\otimes Ext^1(E_{\lambda},E_{\mu})
\end{equation}
We know that $Ext^1(E_{\lambda},E_{\mu})^G=
Hom(E_{\lambda}\otimes\Omega^1, E_{\mu})^G$ is equal
to ${\bf C}$ or to $0$, and when it is equal to ${\bf C}$
then $\mu-\lambda=\xi_j$ for some $j$.
We fix the generator $m_{\mu\lambda}$ of 
$Hom(E_{\lambda}\otimes\Omega^1, E_{\mu})^G$
that takes $v_{\lambda}\otimes n_j$
to $v_{\mu}$, 
indeed  $E_{\lambda}\otimes\Omega^1$
contains a unique summand of multiplicity one isomorphic to 
$E_{\mu}$.
This normalization appears already in \cite{B-K} p. 48.
Hence in order to define an element of
$Hom(V_{\lambda}, V_{\mu})$ $\forall\lambda,\mu$
it is enough to give an element of 
$[E]\in Ext^1(gr E, gr E)^G$ and this is the element
corresponding to $\theta$ of \thmref{starting} (i) according
to the isomorphism of \thmref{homext}.
\end{defn0}

The correspondence $E\mapsto [E]$ is functorial, indeed
a $G$-equivariant map $E\to F$ induces first 
a morphism $gr E\mapsto gr F$
and then a morphism of representations of
${\cal Q}_X$ $[E]\mapsto [F]$.

A direct consequence of \thmref{starting} is

\begin{thm0}
\label{main}
Let $G/P$ be a Hermitian symmetric variety.

(i) For any $G$-homogeneous bundle $E$ we have $m([E])=0$,
where $m$ is the invariant Yoneda morphism recalled in \S 2
$$m\colon Ext^1(gr E, gr E)^G \to Ext^2(gr E, gr E)^G$$

(ii) Conversely for any  $R$-module $F$ and any $e\in
Ext^1(F, F)^G$ such that $m(e)=0$ there exists a $G$-homogeneous bundle $E$
 such
that $gr E=F$ and $e=[E]$.
\end{thm0}

\medskip

{\bf Remark} It is well known, although we do not need it,
 that for any bundle $F$ the usual Yoneda morphism 
$Ext^1(F, F)\to Ext^2(F,F)$ is the 
quadratic part of the Kuranishi morphism. In particular
the invariant Yoneda morphism $Ext^1(F, F)^G\to Ext^2(F,F)^G$
is the invariant piece of the quadratic part of the Kuranishi morphism.
\medskip

{\bf Remark} We recall that the functor $E\mapsto grE$ from
$P-mod$ to $R-mod$ is exact. Our description of the quiver and
\thmref{main}
can be thought roughly as an 
additional structure on $R$-mod that allows one to invert the functor
 $gr$.
 \medskip 
 
 The theorem shows how to define relations in 
${\cal Q}_X$ in order to get an equivalence
of categories. The relations have to reflect the vanishing
$m(e)=0$. 
 We have to remark that since in \defref{defmain}
we have fixed a normalization, the relations in ${\cal Q}_X$
can be changed up to scalar multiplications of the maps involved
(see \corref{commutpn}).

\begin{defn0}
\label{rel12} 
Write $e\in Ext^1(gr E, gr E)^G$ as
$$e=\sum g_{\mu\lambda}m_{\mu\lambda}$$
where
$m_{\mu\lambda}\in Ext^1(E_{\lambda},E_{\mu})^G$ were fixed
in  \defref{defmain}
and $g_{\mu\lambda}\in Hom(V_{\lambda},V_{\mu})$ come from
 the isomorphism (\ref{sumtensor}). The equation $m(e)=0$
 becomes
 $$\sum_{\nu, \lambda} \left(\sum_{\mu}(g_{\nu\mu}g_{\mu\lambda})
(m_{\nu\mu}\wedge m_{\mu\lambda})\right)=0$$
where 
$m_{\nu\mu}\wedge m_{\mu\lambda}\in Ext^2(E_{\lambda},E_{\mu})^G$
is the Yoneda product of $m_{\nu\mu}$, $m_{\mu\lambda}$
and  $g_{\nu\mu}g_{\mu\lambda}\in Hom(V_{\lambda}, V_{\nu})$
are the composition maps. For any fixed $\lambda$ and $\nu$,
the equation
\begin{equation}
\label{relationmg}
\sum_{\mu}(g_{\nu\mu}g_{\mu\lambda})(m_{\nu\mu}\wedge m_{\mu\lambda})=0
\end{equation}
gives a system of at most $\dim Ext^2(E_{\lambda},E_{\mu})^G$
quadratic equations  in the unknowns  $g_{\nu\mu}$ and $g_{\mu\lambda}$

We define the {\bf relations} in ${\cal Q}_X$ as the ideal generated by
all these quadratic equations for any pair $\lambda$ and $\nu$.
\end{defn0}

\begin{thm0}
 \label{teorel}

(i)   For any homogeneous bundle $E$ on $X$  Hermitian 
symmetric variety, $[E]$ 
 satisfies these relations, hence it is a representation of the quiver 
 ${\cal Q}_X$ with relations.

(ii) Conversely given a  representation $e$ of the quiver ${\cal Q}_X$
with relations, there exists a homogeneous bundle $E$ such that
$e=[E]$.
\end{thm0}

{\it Proof}
By definition the relations are equivalent to $\theta\wedge\theta=0$.
 Hence the statement is equivalent
to \thmref{starting} and to \thmref{main} (see also
next \exref{abcd}).
\qedd

\smallskip

The isomorphism class of $[E]$ lives in
$Ext^1(gr E, gr E)^G/Aut^G(gr E)$.
We remark that in each case the isomorphism class of the bundle
determines the isomorphism class of the representation 
of ${\cal Q}_X$ (by the functoriality).
Hence  \thmref{main}
can be reformulated in the following way:

\begin{thm0}
\label{main2}{\bf (Reformulation of \thmref{teorel})}   
Let $X=G/P$ be a Hermitian symmetric variety.
 There is an equivalence of categories among

(i) $G$-homogeneous bundles over $X$. 

(ii) finite dimensional representations of the  quiver (with relations)
 ${\cal Q}_X$ (associating zero to all but a finite number of points of ${\cal Q}_X$).

(iii) Higgs bundles $(F,\theta)$ over $X$.
\end{thm0}

{\bf Subquivers and quotient quivers} Since there is no danger of confusion,
we denote by ${\bf C}{\cal Q}_X$ the path algebra of the quiver with relations ${\cal Q}_X$,
meaning that the algebra has been quotiented by the ideal of relations. 
There are two basic constructions for quiver representations that we will need.

\begin{defn0}
\label{subquot}
Let $gr E=\oplus V_{\lambda}\otimes E_{\lambda}$ so that $V=\oplus V_{\lambda}$ is a ${\bf C}{\cal Q}_X$-module.
For any subspace $V'\subset V$ the submodule generated by $V'$ defines a homogeneous subbundle of $E$.
In case $V'=V_{\lambda'}$ for some $\lambda'$ 
we will call this subbundle {\it the bundle defined by all arrows starting from ${\lambda}'$}.

Also $(V'\colon {\bf C}{\cal Q}_X) :=\{ v\in V| fv\in V'\quad\forall f\in {\bf C}{\cal Q}_X\}$ is a submodule
and the quotient $V/(V'\colon {\bf C}{\cal Q}_X)$ defines a homogeneous quotient of $E$. Let
$\pi_{\lambda'}\colon V\to V_{\lambda'}$ be the projection; in case $V'=Ker~\pi_{\lambda'}$ we have
 $V/(V'\colon {\bf C}{\cal Q}_X)=V/\{ v\in V| \pi_{\lambda'}fv=0\quad\forall f\in {\bf C}{\cal Q}_X\}$
and we will call this quotient bundle {\it the bundle defined by all arrows arriving in $\lambda'$}.
 \end{defn0}

\begin{exa} (compare with \cite{Hi2})
\label{g13}
Let ${\bf P}^3={\bf P}(V)$. The bundle $E=\wedge^2V$ on 
$X=Gr({\bf P}^1,{\bf P}^3)$ has
$gr E={\cal O}(-1)\oplus\Omega^1(1)\oplus{\cal O}(1)$.
The corresponding representation of the quiver
associates to
\[\begin{array}{ccc}
&& {\cal O}(1) \cr
&&\dow{}\cr
{\cal O}(-1)&\lef{}&\Omega^1(1)\cr
\end{array}\]
the diagram of linear maps
\[\begin{array}{ccc}
&&{\bf C}\cr
&&\dow{\theta_1}\cr
{\bf C}&\lef{\theta_2}&{\bf C}\cr
\end{array}\]
Equivalently $\theta$ splits into the two summands
\[\theta_1\colon{\cal O}(1)\otimes\Omega^1\rig{}\Omega^1(1)\]
and
\[\theta_2\colon\Omega^1(1)\otimes\Omega^1\rig{}{\cal O}(-1)\]
and satisfies $\theta\wedge\theta=0$ because
\[Ext^2({\cal O}(1), {\cal O}(-1))^G=
Hom({\cal O}(1)\otimes\Omega^2,{\cal O}(-1) )^G=0\] 
In fact in ${\cal Q}_X$ the commutativity of the diagram
\[\begin{array}{ccc}
0&\lef{}& {\cal O}(1) \cr
\dow{}&&\dow{}\cr
{\cal O}(-1)&\lef{}&\Omega^1(1)\cr
\end{array}\]
is not a relation.
\end{exa}

\begin{thm0}
\label{table}
Let $X$ be an irreducible Hermitian symmetric variety.
The number of connected components of ${\cal Q}_X$
is given by the following table
{\small
\[
\begin{array}{|c|c|c|c|}
\hline
\hbox{Grassmannians}&\hbox{Odd Quadrics}&\hbox{Even Quadrics}
&\hbox{Spinor Varieties}\cr
\hline
SL(n+1)/P(\alpha_{k+1})&Spin(2n+1)/P(\alpha_1)&
Spin(2n+2)/P(\alpha_1)&Spin(2n+2)/P(\alpha_{n+1})\cr
Gr({\bf P}^{k},{\bf P}^n)&Q_{2n-1}\quad\ n\ge 2&Q_{2n}\quad\ n\ge 2&
{\frac{1}{2}}Gr({\bf P}^n,Q_{2n})\quad\ n\ge 3\cr
\hline
n+1&2&4&4\cr
\hline
\end{array}\]
\[\begin{array}{|c|c|c|}
\hline
\hbox{Lagrangian Grassmannians}&\hbox{Cayley Plane}&X_{27}\cr
\hline
Sp(2n)/P(\alpha_n)
&E_6/P(\alpha_1)
&E_7/P(\alpha_1)\cr
Grn({\bf P}^{n-1},{\bf P}^{2n-1})\quad\ n\ge 2&
{\bf OP}^2&\cr
\hline
2&3&2\cr
\hline
\end{array}\]
}
\end{thm0}

{\it Proof} The number of connected components is equal
to the index of the lattice $\langle\xi_1,\ldots ,\xi_m\rangle_{{\bf Z}}$
in $\langle\lambda_1,\ldots, \lambda_n\rangle_{{\bf Z}}$.
It is easy to check in any case that
$$\langle\xi_1,\ldots ,\xi_m\rangle_{{\bf Z}}=\langle\alpha_1,\ldots, \alpha_n\rangle_{{\bf Z}}$$
by the shape of the roots (the list in the exceptional cases
is in \cite{Snow}). Hence the number of connected components
is given in any case by the determinant of the corresponding Cartan 
matrix,
and these are well known (see e.g. \cite{F-H} exerc. 21.18).

\qedd

\medskip

Every homogeneous bundle
$E$ on $X$ splits as $E=\oplus E^{(i)}$
where the sum is over the connected components
of ${\cal Q}_X$, and $gr(E^{(i)})$
contains only irreducible bundles corresponding
to points of the connected component labelled by $i$.
We analyze separately each of the irreducible Hermitian symmetric 
varieties. The decomposition of $E_{\lambda}\otimes\Omega^1$
in the cases where $G$ is of type $A$, $D$ or $E$ appears
already in Prop. 2 of \cite{B-K}.

\noindent$\bullet$ When $G=SL(n+1)$ then $X=G/P(\alpha_{k+1})$
is the Grassmannian $Gr({\bf P}^k, {\bf P}^n)$.
In this case all the roots $\Omega^1_X$ are
$\beta_{ij}=-\sum_{t=i}^j\alpha_t$ for $ 1\le i\le k+1\le j\le n $.
If $U$ and $Q$ are the universal and the quotient bundle, it is well 
known that $\Omega^1=U\otimes Q^*$,
$\Omega^2=[Sym^2U\otimes\wedge^2Q^*]\oplus [\wedge^2U\otimes Sym^2Q^*]$.

Here $\mu(\Omega^1)=-\frac{n+1}{(k+1)(n-k)}$.
Every irreducible bundle on $X$ can be described by
$E={\cal S}^{\alpha}U\otimes{\cal S}^{\beta}Q^*(t)$
for some partitions $\alpha$, $\beta$ and for some $t\in{\bf Z}$.
The $n+1$ connected components are distinguished
by the class of $(|\alpha|, |\beta|)\in {\bf Z}_2\times{\bf Z}_2$
modulo the lattice $\langle (-1,1),(k+1,n-k)\rangle_{\bf Z}$.
If $G.C.D.\left( n+1,(k+1)(n-k)\right)=1$
the components are distinguished more easily by
$(k+1)(n-k)\mu(E)=0,1,\ldots, n (mod\ n+1)$.

\noindent$\bullet$ When $k=0$ we get $X={\bf P}^n$. Due to the importance of this case 
in the applications we stress our attention on it. We saw before  in
\lemref{alphapn} the corresponding roots $\xi_1,\ldots, \xi_n$. 
We have the simple formulas (of course some summands can
 be zero)
\[E_{\lambda}\otimes \Omega^1=\oplus_{i=1}^n E_{\lambda+\xi_i}\]
\[E_{\lambda}\otimes \Omega^2=\oplus_{1\le i<j\le n}
 E_{\lambda+\xi_i+\xi_j}\]

In \corref{commutpn} we will see
 that the relations in the quiver
${\cal Q}_{{\bf P}^n}$ can be summed up by saying that
for any weight $\lambda\in D_1$ and any $1\le i<j\le n$
 all diagrams

\[\begin{array}{ccc}
E_{\lambda+\xi_i}&\lef{}&E_{\lambda}\cr
\dow{}&&\dow{}\cr
E_{\lambda+\xi_i+\xi_j}&\lef{}&E_{\lambda+\xi_j}\cr
\end{array}\]
have to be commutative.
This fits with \cite{B-K}. The quiver ${\cal Q}_{{\bf P}^n}$ is isomorphic to the
half-space of ${\bf Z}^n$ defined by the inequalities
 $x_1\ge x_2\ge\ldots\ge x_n $
for $(x_1,\ldots ,x_n)\in{\bf Z}^n$
with arrows following the standard basis 
(with the directions reversed).
Here $\mu(\Omega^1)=-\frac{n+1}{n}$.
The $n+1$ connected components are distinguished by
$n\mu(E)=0,1,\ldots, n (mod\ n+1)$
for an irreducible $E$.

\noindent$\bullet$ In the case of odd dimensional quadrics 
$Spin(2n+1)/P(\alpha_1)=Q_{2n-1}\subset{\bf P}^{2n}$ we have that 
$\Omega^1$ has maximal weight $-\alpha_1
=2\lambda_2-2\lambda_1$ for $n=2$ and has maximal weight
$-\alpha_1=\lambda_2-2\lambda_1$ for $n\ge 3$, while
$\Omega^2$ has maximal weight $2\lambda_2-3\lambda_1$ for $n=2$,
 $2\lambda_3-3\lambda_1$ for $n=3$ and 
$\lambda_3-3\lambda_1$ for $n\ge 4$. 
Denote again by $\xi_1,\ldots, \xi_m$ ($m=2n-1$) the roots of $\Omega^1$. 
We have
\[E_{\lambda}\otimes \Omega^1=\oplus_{i=1}^m E_{\lambda+\xi_i}\]
while $E_{\lambda}\otimes \Omega^2$
is contained in $\oplus_{1\le i<j\le m}
 E_{\lambda+\xi_i+\xi_j}$ and can be determined according to 
$\lambda$ by the explicit algorithm in \cite{Li}. When $\lambda\gg 0$
then we have the equality.
Here $\mu(\Omega^1)=-1$ and $\mu(S)=-\frac{1}{2}$
for the spinor bundle.
The two connected components are distinguished by
$2\mu(E)=0,1 (mod\  2)$ for an irreducible $E$.

\noindent$\bullet$ In the case of even dimensional quadrics 
$Spin(2n+2)/P(\alpha_1)=Q_{2n}\subset{\bf P}^{2n+1}$ 
($\lambda_n$ and $\lambda_{n+1}$ correspond to the
 two spinor bundles) we have that 
$\Omega^1$ has maximal weight $\lambda_2+\lambda_3-2\lambda_1$
for $n=2$ and $\lambda_2-2\lambda_1$ for $n\ge 3$, while
$\Omega^2$ splits with two maximal weights $2\lambda_2-3\lambda_1$ and
$2\lambda_3-3\lambda_1$ for $n=2$ 
(this is the grassmannian of lines in ${\bf P}^3$ 
already considered), and
it has maximal weight $\lambda_3+\lambda_4-3\lambda_1$
 for $n=3$ and $\lambda_3-3\lambda_1$ for $n\ge 4$. 
Here $\mu(\Omega^1)=-1$ and $\mu(S)=-\frac{1}{2}$
for the two spinor bundles.
The knowledge of $\mu$ is not enough to distinguish
the several components.
If $E=E_{\sum p_i\lambda_i}$
the four components are distinguished by
$[(p_n,p_{n+1})]\in {\bf Z}_2\times{\bf Z}_2$.

\noindent$\bullet$ In the case of spinor variety $Spin(2n+2)/P(\alpha_{n+1})$
we have the universal bundle $U$ of rank $n+1$ and it is well known
that $\Omega^1=\wedge^2 U$ and 
$\Omega^2=\wedge^2(\wedge^2 U)={\cal S}^{2,1,1}U$. Let 
$m={{n+1}\choose 2}$ and let $\xi_1,\ldots ,\xi_m$ be the roots of
 $\Omega^1$.
Then it is easy to check that
\[E_{\lambda}\otimes \Omega^1=\oplus_{i=1}^m E_{\lambda+\xi_i}\]
while $E_{\lambda}\otimes \Omega^2$
is contained in $\oplus_{1\le i<j\le m}
 E_{\lambda+\xi_i+\xi_j}$ that can be determined according to 
$\lambda$ by the classical Littlewood-Richardson rule 
(because the semisimple part of $P(\alpha_{n+1})$ is $SL(n+1)$).
 When $\lambda\gg 0$
then we have the equality.
Here $\mu(\Omega^1)=-\frac{4}{n+1}$.
If G.C.D.(4,n+1)=1 then the four connected components are distinguished by
$(n+1)\mu(E)=0,1,2,3 (mod\  4)$.
Otherwise the knowledge of $\mu$ is not enough to distinguish
the several components.
Every irreducible bundle on $X$ can be described by
$E={\cal S}^{\alpha}U\otimes{\cal O}(t)$
for some partition $\alpha$ and some integer $t$.
The four connected components are distinguished
by the class of $(|\alpha|, t)\in {\bf Z}_2\times{\bf Z}_2$.

\noindent$\bullet$ In the case of lagrangian maximal grassmannians
$Sp(2n)/P(\alpha_n)$ we have the universal bundle $U$ of rank $n$ and it is well known
that $\Omega^1=Sym^2 U$ and 
$\Omega^2=\wedge^2(Sym^2 U)={\cal S}^{3,1}U$. Let 
$m={{n+1}\choose 2}$ and let $\xi_1,\ldots ,\xi_m$ be the roots of 
$\Omega^1$.
In this case
$E_{\lambda}\otimes \Omega^1$ is contained in
$\oplus_{i=1}^m E_{\lambda+\xi_i}$ and the inclusion can be strict.
Indeed also this computation can be done by using the 
 classical Littlewood-Richardson rule. 
Note that we can write the $\xi_i$ as $\gamma_j+\gamma_k$ where $\gamma_j$ are the
weights of $U$. A fortiori 
$E_{\lambda}\otimes \Omega^2$
is contained in $\oplus_{1\le i<j\le n}
 E_{\lambda+\xi_i+\xi_j}$ and it can be determined according to $\lambda$ by the classical Littlewood-Richardson rule.
Here $\mu(\Omega^1)=-\frac{2}{n}$.
The two connected components are distinguished by
$n\mu(E)=0,1 (mod\  2)$ for an irreducible $E$.

\noindent$\bullet$ In the case
of the Cayley plane $E_6/P(\alpha_1)={\bf OP}^2$ (\cite{L-M}, \cite{I-M}) the semisimple
 part of $P(\alpha_1)$
is $Spin(10)$.  $E_{\lambda_2}$ is a twist of
one of the two spinor bundles
and $\Omega^1=E_{\lambda_2}(-2)$.

Hence $\Omega^2=E_{\lambda_3}(-3)$  is irreducible.
Here $\mu(\Omega^1)=-\frac{3}{4}$. 
The three connected components are distinguished by
$4\mu(E)=0,1,2 (mod\  3)$ for an irreducible $E$.
The Cayley plane has an intrinsic interest because it is a Severi 
variety. 

\noindent$\bullet$ Also the $27$-dimensional case
$E_7/P(\alpha_1)$ has $\Omega^1=E_{\lambda_2}(-2)$ 
and $\Omega^2=E_{\lambda_3}(-3)$ both irreducible. 
Here $\mu(\Omega^1)=-\frac{2}{3}$. 
The two connected components are distinguished by
$3\mu(E)=0,1 (mod\  2)$ for an irreducible $E$.

\medskip

The case of the projective plane ${\bf P}^2$ allows an explicit description of some interest.
Let $(x,y)\in{\cal N}\simeq{\bf C}^2$.
Consider the linear maps
given by matrices with coefficients in $\wedge^*{\cal N}$
\[C_k={1\over{k+1}}\cdot \left[\begin{array}{cccc}
x&y&&\cr
&\ddots&\ddots\cr
&&x&y\cr
\end{array}\right] \quad\hbox{\ of size\ }k\times (k+1) \]

\[B_k={1\over k}\cdot \left[\begin{array}{cclc}
-ky\cr
x&-(k-1)y\cr
&\ddots&\ddots\cr
&&(k-1)x&-y\cr
&&&kx\cr
\end{array}\right] \quad\hbox{\ of size\ }(k+1)\times k \]
Now it is easy to check that
\begin{equation}\label{relp2}
C_{k}\wedge C_{k+1}=0\qquad B_{k+1}\wedge B_k=0\qquad
 C_{k+1}\wedge B_{k+1}+B_k\wedge C_k=0
\end{equation}

The interpretation in terms of representations is the following.
The parabolic subgroup $ P(\alpha_1) \subset SL(3)$ has the form
\[ P(\alpha_1)=\left\{\left[\begin{array}{ccc}e&x&y\cr
0&a_{11}&a_{12}\cr
0&a_{21}&a_{22}\end{array}\right]|\quad e\det A=1\right\}\]

The irreducible representation   of $P(\alpha_1)$
corresponding to $Sym^pQ(t)$ is defined by $Sym^pAe^{-t}$.
Consider the derivative 
${\cal P}=Lie P(\alpha_1)\to {\mbox{ gl}}(Sym^p{\bf C}^2)$ 
and call it (with a slight abuse of notation) $Sym^pA-teI$.
The extension $w\in Ext^1(Sym^{k}Q,Sym^{k-1}Q(-1))^G={\bf C}$ defines
a bundle with representation
\begin{equation}
\label{1ext}
\left[\begin{array}{cc}Sym^{k-1}A+eI&wC_k\cr
0&Sym^{k}A\cr
\end{array}\right]
\end{equation}
 where $w$ is a scalar multiple,
and $w=0$  iff the extension splits.

Analogously,  the extension 
$w\in Ext^1(Sym^{k}Q(2),Sym^{k+1}Q)^G={\bf C}$
 defines
a bundle with representation
\begin{equation}
\label{2ext}
\left[\begin{array}{cc}Sym^{k+1}A&wB_k\cr
0&Sym^{k}A-2eI\cr
\end{array}\right]
\end{equation}
where $w$ is a scalar multiple,
which is zero iff the extension splits.
By \thmref{starting} several extensions as in (\ref{1ext}) and
(\ref{2ext}) fit together to give a representation
$\rho$ of ${\cal P}$ if and only if 
 $ \rho_{|N}\wedge  \rho_{|N}=0$ (see the next \exref{abcd}).
We remark that (\ref{relp2}) are equivalent to the fact
that the only relations in ${\cal Q}_{{\bf P}^2}$
are the commutativity ones (see \corref{commutpn}) in all the square diagrams
and the relation $a_2b_1=0$ in the diagrams
\[\begin{array}{ccc}
&& {\cal O}(t)\cr
&&\dow{b_1}\cr
{\cal O}(t-3)&\lef{a_2}&Q(t-2)\cr
\end{array}\]
for any $t\in{\bf Z}$. These last relations
can be seen as the commutativity
in the diagrams
\[\begin{array}{ccc}
0&\lef{}& {\cal O}(t)\cr
\dow{}&&\dow{b_1}\cr
{\cal O}(t-3)&\lef{a_2}&Q(t-2)\cr
\end{array}\]

\begin{exa}
\label{abcd}
We describe explicitly the homogeneous bundle
on ${\bf P}^2={\bf P}(V)$ corresponding to the representation
that associates to 
\[\begin{array}{ccc}
{\cal O} &\lef{}& Q(1)\cr
\dow{}&&\dow{}\cr
Q(-2)&\lef{}&Sym^2Q(-1)\cr
\end{array}\]
the diagram of linear maps
\[\begin{array}{ccc}
{\bf C}^a&\lef{\gamma_1}& {\bf C}^b\cr
\dow{\beta_1}&&\dow{\beta_2}\cr
{\bf C}^c&\lef{\gamma_2}&{\bf C}^d\cr
\end{array}\]
where $a$, $b$, $c$, $d$ are positive integers.
We get
\[
\rho\left[\begin{array}{ccc}e&x&y\cr
0&a_{11}&a_{12}\cr
0&a_{21}&a_{22}\end{array}\right]
=\left[
\begin{array}{cccc}
A^c+2eI&\gamma_2\otimes C_2&\beta_1\otimes B_1&0\cr
0&(Sym^2A)^d+eI&0&\beta_2\otimes B_2\cr
0&0&0&\gamma_1\otimes C_1\cr
0&0&0&A^b-eI\cr
\end{array}\right]\]
and this is a $P$-module
iff (by \thmref{starting}) 
\[\left[\begin{array}{cccc}
0&\gamma_2\otimes C_2&\beta_1\otimes B_1&0\cr
0&0&0&\beta_2\otimes B_2\cr
0&0&0&\gamma_1\otimes C_1\cr
0&0&0&0\cr
\end{array}\right]\wedge
\left[\begin{array}{cccc}
0&\gamma_2\otimes C_2&\beta_1\otimes B_1&0\cr
0&0&0&\beta_2\otimes B_2\cr
0&0&0&\gamma_1\otimes C_1\cr
0&0&0&0\cr
\end{array}\right]=0\]
which is equivalent by (\ref{relp2}) to
\[\gamma_2\cdot\beta_2-\beta_1\cdot\gamma_1=0\]
confirming the commutativity relations.
In the special case $a=b=c=d=1$
and all the maps given by the identity this bundle is $ad V$. 

\end{exa}

The isomorphism classes of representations are equivalent to the orbits in $m^{-1}(0)$ with respect to the
$Aut_G(gr E)$-action.

\section{Computation of Cohomology}

In all this section $X$ is a Hermitian symmetric variety of ADE type.

We want to describe now how to compute the cohomology of a homogeneous bundle $E$ on $X$ from the 
representation of the quiver.

We need the following easy lemma.

\begin{lemma0}{\cite{C-E} lemma XV 1.1}
\label{ce}

 Let the following diagram be commutative
$$\begin{array}{ccccc}
&&C\\
&\near{}&\dow{\phi}&\sear{\psi}\\
A'&\rig{\phi'}&A&\rig{\eta}&A''\\
\end{array}$$ 
 
 and let the row be exact.
 Then $$Im\ \phi /Im \phi' \simeq Im\ \psi$$
\end{lemma0}

Let $$0=E_0\subset E_1\subset E_2\subset\ldots\subset E_r=E$$
be a filtration of a vector bundle (not necessarily homogeneous).

Let now

$$Z^p_j:=Ker\left(H^{j}(E_{p+1}/E_p)\rig{\partial}H^{j+1}(E_{p})\right)$$

$$B^p_j:=Im\left(H^{j-1}(E/E_{p+1})\rig{\partial}H^j(E_{p+1}/E_{p})\right)$$

where the maps are the boundary maps of the two obvious exact sequences.

The following proposition follows from the discussion at the beginning of chapter XV of
\cite{C-E}. For the convenience of the reader we sketch the proof.

\begin{thm0}
\label{filt}
$B^p_j\subset Z^p_j$ and
$$H^j(E)\simeq \oplus_{p=0}^{r-1} Z^p_j/B^p_j $$
\end{thm0}

{\bf Proof}
We have the commutative diagram
 $$\begin{array}{ccccc}
&&H^j(E_{p+1})\\
&\nearrow&\dow{\phi}&\searrow{\psi}\\
H^{j-1}(E/E_{p+1})&\rig{\phi'}&H^j(E_{p+1}/E_{p})&\rig{\eta}&H^{j}(E/E_{p})\\
\end{array}$$

hence $B^p_j\subset Z^p_j$
and from \lemref{ce} we get

\begin{equation}
\label{elena1}
Im\left(H^j(E_{p+1})\rig{\psi}H^{j}(E/E_{p})\right)\simeq Im~\phi/Im~\phi'=
Im~\phi/Ker~\eta=Z^p_j/B^p_j
\end{equation}

Consider also the diagram
 $$\begin{array}{ccccc}
&&H^j(E_{p+1})\\
&\near{}&\dow{\phi_p}&\sear{\psi}\\
H^{j}(E_{p})&\rig{{\phi}_{p-1}}&H^j(E)&\rig{\eta_p}&H^{j}(E/E_{p})\\
\end{array}$$

we get  again from \lemref{ce}

\begin{equation}
\label{elena2}
Im\left(H^j(E_{p+1})\rig{\psi}H^{j}(E/E_{p})\right)\simeq Im({\phi}_p)/Im(\phi_{p-1})
\end{equation}
and since we have  the graduation

$$H^j(E)\simeq \oplus_p Im({\phi}_p)/Im({\phi}_{p-1})\stackrel{(\ref{elena1})(\ref{elena2})}{=}\oplus_p Z^p_j/B^p_j$$
we get the result.\qedd  

\bigskip

We return now to the case of homogeneous bundles.

We need a short digression about homogeneous bundles whose quiver representation has support
on an $A_n$-type quiver, that is $gr E=\oplus V_\lambda\otimes E_{\lambda}$ and $V_{\lambda}$
is zero outside a path connecting the vertices $\{\lambda+p\xi_j | 0\le p\le k\}$.

The following theorem is well known since the former work on quivers by P. Gabriel (see \cite{G-R}).

\begin{thm0}
\label{segrehomog}
 Every representation of the $A_m$-quiver is the direct sum
of irreducible representations with dimension vector
$$(0,0,\ldots,0,1,1,\ldots, 1,0,\ldots,0)$$
where the nontrivial linear maps are isomorphisms.
\end{thm0}

The reader can enjoy to deduce the previous theorem as a consequence of \thmref{main2} for $X={\bf P}^1$
and the Segre-Grothendieck theorem,
which says that every bundle on ${\bf P}^1$ splits as the sum of line bundles.

\begin{prop0}
\label{needADE}
Let $E_{\lambda}$ and $E_{\mu}$ be in two adjacent Bott chambers with
$H^{i}(E_{\lambda})\simeq H^{i+1}(E_{\mu})\simeq W$, then
$\mu-\lambda=k\xi_j$ for some integer $k$ and some root $\xi_j$ of $\Omega^1$.
We have

$$\dim Hom\left(E_{\lambda}\otimes Sym^k\Omega^1,E_{\mu}\right)^G=1$$
\end{prop0}

{\bf Proof}
By (\ref{littelmann}) it is enough to show that there are no other weights 
among $\{a_1\xi_{i_1}+\ldots +a_h\xi_{i_h}| \sum a_i=k\}$
which are equal to $\mu-\lambda$.
With the ADE assumption, $\xi_j$ is a vertex of the convex polytope containing the weights of 
$\Omega^1$, because all the roots have the same length.
Hence $k\xi_j$ is a vertex of the convex polytope containing the weights of 
$Sym^k\Omega^1$.  \qedd

\bigskip

\begin{prop0}
\label{needADE2}
Let $\xi_j$ be a weight of $\Omega^1$.
$$Ext^2(E_{\lambda}, E_{\lambda+2\xi_j})^G=Hom((E_{\lambda}\otimes\Omega^2, E_{\lambda+2\xi_j})^G=0$$
\end{prop0}
{\bf Proof}
Since  $\xi_j$ is a vertex of the convex polytope containing the weights of 
$\Omega^1$, there are no distinct weights $\xi_p$, $\xi_q$ of $\Omega^1$
such that $\xi_j=\frac{1}{2}(\xi_p+\xi_q)$.
Then apply \thmref{homext}. \qedd 

\bigskip

{\bf Remark} Without the ADE assumption the above two propositions are false. For example
if $X=Q_3$, the weights of $\Omega^1$ 
are $\alpha_1$, $\alpha_1+\alpha_2$, $\alpha_1+2\alpha_2$. $\alpha_1+\alpha_2$ is shorter,
indeed $2(\alpha_1+\alpha_2)$ coincides the sum of the two vertices  $(\alpha_1)+(\alpha_1+2\alpha_2)$.
In particular $\dim Hom(Sym^2\Omega^1, E_{2(\alpha_1+\alpha_2)})^G=2$
and $\dim Ext^2({\cal O}, E_{2(\alpha_1+\alpha_2)})^G=1$.
Hence there is no indecomposable homogeneous bundle $E$
with support $A_2$ such that $gr E=\oplus_{i=0}^2E_{i(\alpha_1+\alpha_2)}$.

\bigskip

With the assumption of the two previous propositions, note that the distinguished elements in
$Hom\left(E_{\lambda+p\xi_j}\otimes \Omega^1,E_{\lambda+(p+1)\xi_j}\right)^G$
which were chosen in \defref{defmain} give a distinguished element in
$Hom\left(E_{\lambda}\otimes Sym^k\Omega^1,E_{\mu}\right)^G$,
which is one dimensional by \propref{needADE}.
These elements allow one to define extensions of the form
$$0\rig{}E_{\lambda+(p+1)\xi_j}\rig{}Z_p\rig{}E_{\lambda+p\xi_j}\rig{}0$$
which fit together (by \thmref{main} (ii), since the corresponding $Ext^2$ vanish by \propref{needADE2})
 giving a bundle $P'$
with $gr E=\oplus_{p=0}^{k-1} E_{\lambda+p\xi_j}$
and two exact sequences
(this argument is similar to the one in \cite{Dem})
\begin{equation}
\label{zp1}
0\rig{}Z'\rig{}P'\rig{} E_{\lambda}\rig{}0
\end{equation}

\begin{equation}
\label{zp2}
0\rig{} E_{\mu}\rig{}Z'\rig{}Z'/ E_{\mu}\rig{}0
\end{equation}

\begin{thm0}
\label{spaghettinat}

$$H^j(P')=0\quad\forall j$$
\end{thm0}

We need a short preparation in order to prove \thmref{spaghettinat}.
Let $\lambda'$ (resp. $\mu'$) be the vertex of the Bott chamber containing $\lambda$ (resp. $\mu$).
Let $A$ be the unique indecomposable bundle in the extension

$$0\rig{}E_{\mu'}\rig{}A\rig{}E_{\lambda'}\rig{}0$$

\begin{prop0}
\label{azero}
$$H^i(A)=0\quad\forall i$$
\end{prop0}

{\bf Proof}
The boundary map 
$H^i(E_{\lambda'})\rig{}H^{i+1}(E_{\mu'})$
can be seen as the cup product of class of the Schubert cell
corresponding to $E_{\lambda'}$ as subbundle of $\Omega^i$ (by Hodge theory)
with the hyperplane class  in $H^1(\Omega^1)$
and it is nonzero by
\cite{Hi} Coroll. V 3.2.    \qedd

\medskip

\begin{prop0}
\label{only}
$gr (E_{\mu'}\otimes W)$ contains only $E_{\mu}$ as direct summand with
$H^*\simeq W$.
\end{prop0}

{\bf Proof} Let $E_{\alpha}$ be the irreducible bundle such that $H^0(E_{\alpha})=W$.
The weights of $W$ as $G$-module lie in a convex polytope $P_W$ whose vertices are the reflections
of $\alpha$ through the hyperplanes $H_{\phi}$ (for any root $\phi$ of $G$) 
which separate the Weyl chambers of $G$ (see\cite{F-H} pag. 204).
The weights of $E_{\mu'}\otimes W $ lie inside $P_W+\mu'$.

Let $\tilde P_W$ be the convex polytope whose vertices are the reflections
of $g+\alpha$ through the hyperplanes $H_{\phi}$. Note that $P_W$ is strictly
contained in $\tilde P_W$ and there is a natural bijective correspondence $f$ between
the vertices of $P_W$  and the vertices of $\tilde P_W$ such that if $\beta\gamma$
is an edge of $P_W$ of length $d\sqrt{2}$ then $f(\beta)f(\gamma)$
is a parallel edge of $\tilde P_W$ of length $(d+1)\sqrt{2}$.
Precisely the corresponding vertices $ \tilde \beta $ and $ \beta$  
respectively of $\tilde P_W$ and $P_W$ differ by $w_{\beta}(g)$
for a composition of reflections  $w_{\beta}$  defined by  $\beta=w_{\beta}(\alpha)$.
The point of $P_W$ of least distance from $\tilde \beta $ is $\beta$. 

We have that 
$\mu'=w(g)-g$ for some $w$. Let $\overline{\mu}=w(\alpha)$, then $w=w_{\overline{\mu}}$.
 Then $\mu=w_{\overline{\mu}}(\alpha+g)-g=\overline{\mu}+\mu'$ is a vertex of  $P_W+\mu'$ ,
 hence it is a maximal weight
of  $E_{\mu'}\otimes W $.

By Bott theorem all the weights ${\nu}$ such that $H^i(E_{\nu})=W$ for some $i$ are obtained from
 $\alpha$ after reflecting through the hyperplanes which separate the Bott chambers of $G$.
All these weights are some of the vertices of $\tilde P_W-g$.

It is enough to show that the vertices of $\tilde P_W-g$ meet $P_W+\mu'$ only in the point $\mu$.

The distance 
of $\tilde\beta -g$ from $P_W-g+(\mu'+g)=P_W-g+w_{\overline{\mu} }(g)$ vanishes only when 
$\tilde\beta -w_{\overline{\mu}}(g)\in P_W$ and
this happens iff $w_{\beta}(g) = w_{\overline{\mu}}(g) $ (since the point of least distance between $\tilde \beta$
 and $  P_W$ is $\tilde\beta-w_{\beta}(g)$);
thus $\beta =\overline{\mu}$. Then ${\tilde \beta}-g=\beta+w_{\beta}(g)-g=\beta+w_{\overline{\mu}}(g)-g=
\overline{\mu}+\mu'=\mu$.
\qedd

\medskip

{\bf Proof of \thmref{spaghettinat}}
Let now $K$ be the submodule in $A\otimes W$ generated by the direct summands isomorphic to $E_{\lambda}$
(it can be shown that there is only one but we do not need this fact).
We have the exact sequence

$$0\rig{}K\rig{}A\otimes W\rig{}Q\rig{}0$$

By \propref{only} we have that $H^j(K)^W$ and $H^j(Q)^W$ are nonzero at most for $j=i$
or $j=i+1$. 

We claim that $gr K$ contains all the direct summands isomorphic to $E_{\mu}$, otherwise $E_{\mu}\subset gr Q$, 
and we would have
$H^{i+1}(Q)^W\neq 0$, hence by \propref{azero}
$H^{i+2}(K)^W\neq 0$ which is a contradiction.
Hence we get $H^j(Q)^W=0\quad\forall j$ and it follows
$$H^j(K)^W=0\quad\forall j$$

At last, let $S'$ be the quotient of $K$ obtained restricting the quiver representation
to the path joining the vertices corresponding to $E_{\lambda}$ and $E_{\mu}$.

We have 
$$0\rig{}K'\rig{}K\rig{}S'\rig{}0$$

Now $H^j(gr K')^W=0\quad\forall j$, hence $H^j(K')^W=0\quad\forall j$
and it follows
$H^j(S')=0\quad\forall j$.
Decompose $S'$ into its irreducible components (see \thmref{segrehomog}), we get that $S'$ is
isomorphic to the direct sum of several copies of $P'$, by the definition of $K$. \qedd

\medskip

From the sequence (\ref{zp1}) and \thmref{spaghettinat} we have the isomorphism
$$H^j(E_{\lambda})\rig{\partial}H^{j+1}(Z')^W$$
and from (\ref{zp2}) an isomorphism
$$H^{j+1}(E_{\mu})\rig{\simeq}H^{j+1}(Z')$$
hence we get a distinguished isomorphism
\begin{equation}
\label{disting}
j_{\mu\lambda}\colon H^j(E_{\lambda})\rig{}H^{j+1}(E_{\mu})
\end{equation}
  
\begin{lemma0}
\label{spagmagro}
Let $E_{\lambda}$ and $E_{\mu}$ be in two adjacent Bott chambers with
$H^{j-1}(E_{\lambda})\simeq H^j(E_{\mu})\simeq W$.
Denote by $P$ the homogeneous bundle corresponding to the $A_n$-type, starting from $E_{\lambda}$
and arriving in $E_{\mu}$, with the same representation quiver maps as for $E$
(it exists by \thmref{main} (ii) again, by the same argument as before).  Then  the boundary map
$$W\otimes V_{\lambda}=H^{j-1}(P/V_{\mu}E_{\mu})^W\rig{\partial} H^j(V_{\mu}E_{\mu})=W\otimes  V_{\mu} $$
is the tensor product of the distinguished isomorphism in (\ref{disting}) and
the composition of the maps of the quiver representation.
\end{lemma0}
{\bf Proof}
We first prove the theorem for $P$ irreducible. We may assume $\dim V_{\lambda+p\xi_j}=1$
for  $0\le p\le k$ and $\lambda+k\xi_j =\mu$, moreover
$P$ defines nonzero elements in the one dimensional spaces
$$Hom\left(V_{\lambda+p\xi_j}\otimes E_{\lambda+p\xi_j}\otimes\Omega^1,
V_{\lambda+(p+1)\xi_j}\otimes E_{\lambda+(p+1)\xi_j}\right)^G=$$
$$=Hom(V_{\lambda+p\xi_j}, V_{\lambda+(p+1)\xi_j})\otimes 
Hom\left(E_{\lambda+p\xi_j}\otimes\Omega^1,E_{\lambda+(p+1)\xi_j}\right)^G$$
There is a natural isomorphism between $$\bigotimes_{i=0}^{k-1}
Hom\left(V_{\lambda+p\xi_j}\otimes E_{\lambda+p\xi_j}\otimes\Omega^1,
V_{\lambda+(p+1)\xi_j}\otimes E_{\lambda+(p+1)\xi_j}\right)^G$$
and
$$Hom\left(V_{\lambda}\otimes E_{\lambda}\otimes Sym^k\Omega^1,
V_{\mu}\otimes E_{\mu}\right)^G=
Hom(V_{\lambda}, V_{\mu})\otimes 
Hom\left(E_{\lambda}\otimes Sym^k\Omega^1,E_{\mu}\right)^G$$

where in  $Hom(V_{\lambda}, V_{\mu})$ we perform the composition of the quiver representation maps.

It is clear that the element obtained in 
$Hom(V_{\lambda}, V_{\mu})\otimes 
Hom\left(E_{\lambda}\otimes Sym^k\Omega^1,E_{\mu}\right)^G$
 is enough to reconstruct $P$.
 
Now we consider the two exact sequences
$$0\rig{}Z\rig{}P\rig{}V_{\lambda}\otimes E_{\lambda}\rig{}0$$
$$0\rig{}V_{\mu}\otimes E_{\mu}\rig{}Z\rig{}P'\rig{}0$$
From the first sequence we have $$H^j(E_{\lambda}\otimes V_{\lambda})\rig{\partial}H^{j+1}(Z)^W$$
and from the second one an isomorphism (by \thmref{spaghettinat})
$$H^{j+1}(E_{\mu}\otimes V_{\mu})\rig{\simeq}H^{j+1}(Z)$$
hence we get a map
\begin{equation}
\label{disting2}
c_{\mu\lambda}\colon H^j(E_{\lambda}\otimes V_{\lambda})\rig{}H^{j+1}(E_{\mu}\otimes V_{\mu})
\end{equation}
which by the construction is the tensor product of the distinguished isomorphism $j_{\mu\lambda}$ constructed
 in (\ref{disting}) and the composition  of the maps of the quiver representation, as we wanted.
 
 In general we have $P=\oplus P_i$ where $P_i$ are irreducible by \thmref{segrehomog}.
Moreover we have $V_{\lambda}=\oplus V_{\lambda}^i$, $V_{\mu}=\oplus V_{\mu}^i$
where every $V_{\lambda}^i$ and $V_{\mu}^i$ has dimension one or zero and for each $i$
the morphism $W\otimes V_{\lambda}^i=H^{j-1}(P_i/V_{\mu}^iE_{\mu})^W\rig{\partial} H^j(V_{\mu}^iE_{\mu})=
W\otimes V_{\mu}^i $
coincides again with the tensor product of the distinguished isomorphism $j_{\mu\lambda}$ constructed
 in (\ref{disting}) and the composition  of the maps of the quiver representation. \qedd
\bigskip

We construct now maps $H^j(gr E)\rig{c_j}H^{j+1}(gr E)$
by patching together the maps
$c_{\mu\lambda}$
already constructed in (\ref{disting2}), that is
\begin{defn0}
\label{defc}
$$c_j:=\sum c_{\mu\lambda}$$
where the sum is all over pairs $\lambda$, $\mu$  in two adjacent Bott chambers and
$H^{j}(E_{\lambda})\simeq H^{j+1}(E_{\mu})$.
\end{defn0}
Although separately the isomorphism $j_{\mu\lambda}$ in (\ref{disting})
 and the composition of the quiver representations maps
depend on the choices made in \defref{defmain},  it is easy to check that their tensor product
does not depend on these choices (the scalar multiple that one has to change cancel together). 
Moreover the construction in \defref{defc} is functorial, that is given $E\rig{}F$
we get a morphism 
$H^*(gr E)\rig{}H^*(gr F)$. We see now that $H^*(gr E)$ is a complex and it
gives a way to compute the cohomology.

\smallskip

In the case of ${\bf P}^n$ this construction can be made more explicit.
We have maps given by
$g_{\lambda,i}\colon W_{\lambda}\to W_{\lambda+\alpha_1+\ldots +\alpha_{i+1}}=
W_{\lambda'}$.
Let $\lambda=\sum_{i=1}^np_i\lambda_i$.
Let $p_i({\lambda})=-\sum_{j=1}^{i+1}(p_j+1)$.
Composing the  maps 
$W_{\lambda+j(\alpha_1+\ldots +\alpha_{i+1})}\to
W_{\lambda+(j+1)(\alpha_1+\ldots +\alpha_{i+1})}$
for $i$ fixed and $j=0,\ldots,p_i-1$ we get
$W_{\lambda}\to W_{\lambda+p_i(\lambda) (\alpha_1+\ldots +\alpha_{i+1}) }$
 and we get
$W_{\lambda}\rig{g'_{\lambda,i}} W_{\lambda '}$ where
$H^i(E_{\lambda})=H^{i+1}(E_{\lambda'})$ and
$g'_{\lambda,i}=\prod_{j=1}^{p_i(\lambda)}g_{\lambda+(j-1)
(\alpha_1+\ldots +\alpha_{i+1}),i }$. 
The corresponding maps $c_0$, $c_1$ in the case of ${\bf P}^2$
are shown in the following picture.

\hspace*{2cm}
\includegraphics[scale=0.3]{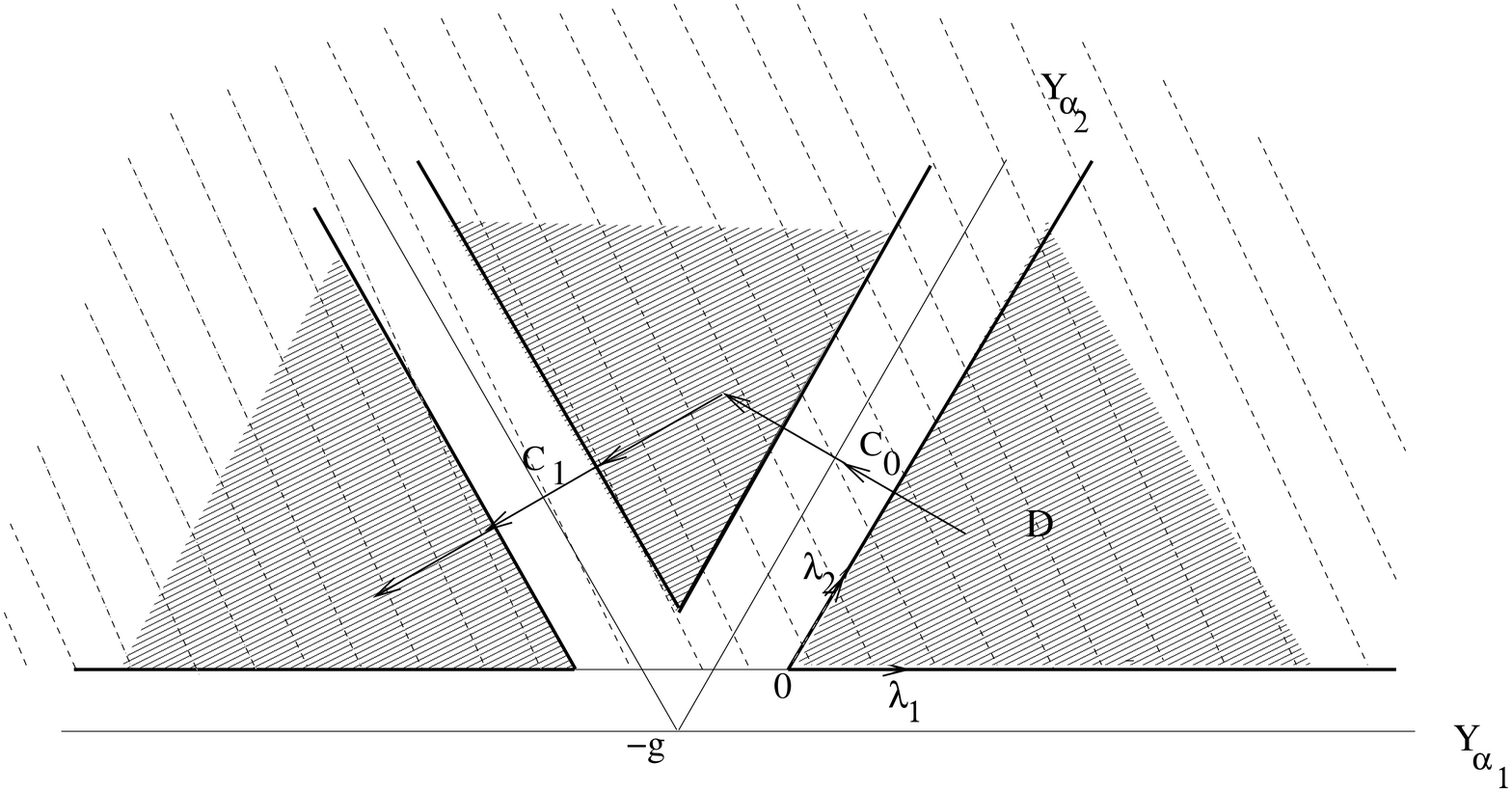}

\bigskip

{\bf Remark} We warn the reader that the use of the distinguished isomorphism (\ref{disting})
is not a formal and superfluous addition, but it determines the correct signs
which are necessary in concrete computations. For example, assume we have
$\lambda$, $\mu$, $\nu$ in three consecutive adjacent Bott chambers such that
$H^{j}(E_{\lambda})\simeq H^{j+1}(E_{\mu})\simeq H^{j+2}(E_{\nu})$,
and $\lambda$, $\mu'$, $\nu$ in the same situation (at most two $\mu$'s exist between
$\lambda$ and $\nu$), it may be shown as an application
of the well known relation in \cite{C-E} III prop. 4.1 that we have the anticommutativity relation
$$j_{\nu\mu}j_{\mu\lambda}=-j_{\nu\mu'}j_{\mu'\lambda}$$
The next theorem implies in this case that
$$c_{\nu\mu}c_{\mu\lambda}=-c_{\nu\mu'}c_{\mu'\lambda}$$
hence it follows by the construction of $c_{\mu\lambda}$ that the corresponding composition of quiver representation
maps is commutative for the square
$$\begin{array}{ccc}
\lambda&\rig{}&\mu\\
\dow{}&&\dow{}\\
\mu'&\rig{}&\nu\\
\end{array}
$$

taken from the Hasse quiver.
In the last section about Olver maps we give more informations in the case of Grassmannians. 

\begin{thm0}
\label{filt1}
$\left( H^*(gr E), c_*\right)$ is a complex and its cohomology
is given by $$\frac{Ker~c_i}{Im~c_{i-1}}=H^i(E)$$
\end{thm0}

{\bf Proof}
Let $W$ be any irreducible $G$-module and let $n=\dim X$.
It is enough to compute that 
$$H^j(E)^W=
\frac{Ker\left(H^i(gr E)^W\rig{c_i}H^{i+1}(gr E)^W\right)}
{Im\left(H^{i-1}(gr E)^W\rig{c_{i-1}}H^{i}(gr E)^W\right)}
$$

We consider the filtration of $E$ defined in the following way.

$E_1$ is defined by taking all arrows starting from any $F\in gr E$ such that $H^n(F)^W\neq 0$
(see \defref{subquot}). 

$E_2$ is defined taking all arrows starting from any $F\in gr E$ such that $$H^n(F)^W\oplus H^{n-1}(F)^W\neq 0$$

In general $E_i$ is defined taking all arrows starting from any $F\in gr E$ such that 
$$\oplus_{j=0}^{i-1} H^{n-j}(F)^W\neq 0$$

We get
$$H^j(gr E_{i+1}/E_i)^W=\left\{\begin{array}{ll}
H^{n-i}(gr E)^W&\hbox{if\ }j=n-i\\
0&\hbox{if\ }j\neq n-i\\
\end{array}\right.$$
hence by the spectral sequence
$$H^j(E_{i+1}/E_i)^W=\left\{\begin{array}{ll}
H^{n-i}(gr E)^W&\hbox{if\ }j=n-i\\
0&\hbox{if\ }j\neq n-i\\
\end{array}\right.$$

We have the commutative diagram

$$
\begin{array}{ccccc}
H^{i-1}(E_{n-i+2}/E_{n-i+1})^W\\
||\\
H^{i-1}(gr E/E_{n-i+1})^W\\
\dow{f}&\sear{\partial}\\
H^{i-1}(E/E_{n-i+1})^W&\rig{\partial}&H^i(E_{n-i+1}/E_{n-i})^W
&\rig{\partial}&H^{i+1}(E_{n-i})^W\\
&&&\sear{\partial}&\dow{g}\\
&&&&H^{i+1}(gr E_{n-i})^W\\
&&&&||\\
&&&&H^{i+1}(E_{n-i}/E_{n-i-1})^W\\
\end{array}$$

where $f$ is the projection given by the spectral sequence ($H^i(gr E/E_{n-i+1})^W=0$) and
$g$ is injective (because $H^i(gr E_{n-i})^W=0$). Moreover we remark that the central term is  
$$H^i(E_{n-i+1}/E_{n-i})^W=H^i(gr E)^W$$

It follows from this diagram and \thmref{filt} that

$$H^i(E)^W=Z^{n-i}_i/B^{n-i}_i=\frac{Ker
\left( H^{i}(E_{n-i+1}/E_{n-i})^W\rig{\partial} H^{i+1}(E_{n-i}/E_{n-i-1})^W \right)}
{Im \left(  H^{i-1}(E_{n-i+2}/E_{n-i+1})^W\rig{\partial} H^{i}(E_{n-i+1}/E_{n-i})^W \right)}
$$

Now it is enough to show that the boundary map
$$ H^{i-1}(E_{n-i+2}/E_{n-i+1})^W\rig{\partial} H^{i}(E_{n-i+1}/E_{n-i})^W $$
induced by the exact sequence
$$0\rig{}E_{n-i+1}/E_{n-i}\rig{}E_{n-i+2}/E_{n-i} \rig{}E_{n-i+2}/E_{n-i+1}\rig{}0$$
is the composition of the quiver representation maps tensored with $j_{\mu\lambda}$
in (\ref{disting}).

\lemref{spagmagro} tells that this is true in the particular case
of quiver representations with support $A_m$, and we will bring back to that case.
Pick $V_{\lambda}\otimes E_{\lambda}\subset gr E_{n-i+2}/E_{n-i+1}$
and
$V_{\mu}\otimes gr E_{\mu}\subset E_{n-i+1}/E_{n-i}$
such that $W\simeq H^{i-1}(E_{\lambda})\simeq H^{i}(E_{\mu})$.

We have to show that the composition

$$H^{i-1}(V_{\lambda}\otimes E_{\lambda})\rig{l}H^{i-1}(E_{n-i+2}/E_{n-i+1})^W\rig{\partial} H^{i}(E_{n-i+1}/E_{n-i})^W
\rig{} H^{i}(V_{\mu}\otimes E_{\mu})$$

is obtained by composing the maps appearing in the quiver representation
from $V_{\lambda}$ to $V_{\mu}$.

Consider the commutative diagram

$$\begin{array}{ccccccccc}
&&0&&0&&0\\
&&\dow{}&&\dow{}&&\dow{}\\
0&\rig{}& K\cap (E_{n-i+1}/E_{n-i}) &\rig{} & E_{n-i+1}/E_{n-i}&\rig{}&Q'&\rig{}&0\\
&&\dow{}&&\dow{}&&\dow{}\\
0&\rig{}&K&\rig{} &E_{n-i+2}/E_{n-i}&\rig{}&Q&\rig{}&0\\
&&&&\dow{}&&\dow{}\\
&&&&E_{n-i+2}/E_{n-i+1}&\rig{}&Q''\\
&&&&\dow{}&&\dow{}\\
&&&&0&&0\\
\end{array}$$

where $Q$ is the quotient of $E_{n-i+2}/E_{n-i}$ 
obtained by taking all arrows arriving in $E_{\mu}$ (see \defref{subquot})
and the other bundles are defined
from the diagram itself.

This diagram induces the diagram

$$\begin{array}{ccccccc}
&&H^{i}(V_{\lambda}\otimes E_{\lambda})\\
&&\dow{}&\sear{}\\
&&H^{i}(E_{n-i+2}/E_{n-i+1})^W &\rig{}&H^{i}(Q'')^W \\
&&\dow{\partial}&&\dow{\partial}\\
H^{i+1}(K\cap (E_{n-i+1}/E_{n-i}) )^W&\rig{f}&H^{i+1}
(E_{n-i+1}/E_{n-i})^W& \rig{} &H^{i+1}(Q')^W\\
&&\dow{h}\\
&&H^{i+1}(V_{\mu}\otimes E_{\mu})\\
\end{array}$$

The composition  $hf$ is zero because $E_{\mu}$ is not a vertex of $K$,
then the map $h$ lifts to

$$\begin{array}{ccccccc}
H^{i}(V_{\lambda}\otimes E_{\lambda})\\
\dow{}&\sear{r}\\
H^{i}(E_{n-i+2}/E_{n-i+1})^W &\rig{}&H^{i}(Q'')^W &=&H^{i}(V_{\lambda}\otimes E_{\lambda}) \\
\\
\dow{\partial}&&\dow{\partial}\\
H^{i+1}(E_{n-i+1}/E_{n-i})^W& \rig{} &H^{i+1}(Q')^W\\
\dow{h}&\swar{g}\\
H^{i+1}(V_{\mu}\otimes E_{\mu})\\
\end{array}$$

The last step is to construct the subbundle $P$ of $Q$ taking all arrows
starting from $\lambda$ (see \defref{subquot}), hence $P$ is as in the assumptions of \lemref{spagmagro}.
We get the commutative diagram
$$\begin{array}{ccc}
H^{i}(P/P\cap Q')^W\simeq V_{\lambda}\otimes W&\rig{r}&H^{i}(Q'')^W  \\
\dow{\partial}&&\dow{\partial}\\
H^{i+1}(V_{\mu}\otimes E_{\mu})^W&\rig{k}&H^{i+1}(Q')^W\\
&&\dow{g}\\
&&H^{i+1}(V_{\mu}\otimes E_{\mu})\\
\end{array}$$
where $k$ and $r$ are induced by the inclusions. By the construction of $Q$ we have
$H^{i+1}( gr Q')^W=H^{i+1}(V_{\mu}\otimes E_{\mu})$, 
hence it follows $H^{i+1}(Q')^W=H^{i+1}(V_{\mu}\otimes E_{\mu})$
where the equality is given by $g$ and the composition $gk$ is the identity.

By \lemref{spagmagro} the map ${\partial}$ in the first column of the last diagram is the composition
of the quiver representation maps tensored with $j_{\mu\lambda}$
in (\ref{disting}), then by chasing in the two above diagrams the claim is proved and the proof 
is complete.
\qedd
\bigskip

{\bf Remark} The fact that $(H^*(gr E),c)$ is a complex should be in principle a consequence of the relation
$\theta\wedge\theta=0$. 
Conversely \thmref{filt1} shows that the relation $\theta\wedge\theta=0$, which is quite difficult to be handed directly,
has simpler consequences.
The reader will find some informations more on this topic in the last section about Olver maps.

\bigskip

{\bf Remark} The computation of cohomology allows a interpretation involving the Hasse quiver ${\cal H}_X$
(see section 2). ${\cal H}_X$ is obviously levelled according to \defref{quiverdef}.
Let $\lambda\rig{}\mu\rig{}\nu$ any composition of arrows in  ${\cal H}_X$.
We define quadratic relations in ${\cal H}_X$
asking that the sum of all the composition of two arrows between $\lambda$ and $\nu$ is zero, for all 
$\lambda$ and $\nu$.  Now given a homogeneous bundle $E$ and a irreducible $G$-module $W$
we define a representation of ${\cal H}_X$ in the following way. Let $gr E=\oplus V_{\lambda}E_{\lambda}$.
Given the vertex $\mu$ in ${\cal H}_X$, there is a unique $\lambda$ in the Bott chamber
with vertex $\mu$ such that $H^*(E_{\lambda})\simeq W$. Then we associate to this vertex the $G$-module
$W\otimes V_{\lambda}$. The maps $c_i$ of the complex $H^*(gr E)$ give the maps of this representation.
The direct sum of all these representations for any irreducible $G$-module give a representation of ${\cal H}_X$,
which satisfies the relations we have defined just because  $H^*(gr E)$ is a complex.

So we have constructed a functor from representations of ${\cal Q}_X$ (in finite dimensional vector spaces) 
to representations of ${\cal H}_X$ (in finite dimensional $G$-modules). This functor is not
injective on the objects because the singular weights give zero contribution.
It is easy to see that this functor is neither surjective, so that the representations which are in the image
of the functor make an interesting subcategory.

\medskip

We have  that
for any homogeneous bundle $E$ (on $X$ Hermitian symmetric variety)
the Yoneda product with $[E]\in Ext^1(gr E, gr E)^G$ defines a complex
\[\ldots\rig{}H^i(gr E)\rig{c_i[1]} H^{i+1}(gr E)\rig{}\ldots\]
It is a complex because $m([E])=0$.
We get a functor
from $P$-$mod$ to the (abelian) category $Kom(G$-$mod)$ of complexes
of $G$-modules  $$E\mapsto H^*(gr E)$$ 
It is straightforward to check, by using the properties
of the Yoneda product, that it is an exact functor.
So it is natural to ask about the cohomology of the above complex.
It turns out that, in the ADE case, it gives only the first step of a filtration
of the cohomology $H^*(E)$.
In fact, for any integer $n$, we can consider
 the map $H^i(gr E)\rig{c_i[n]} H^{i+1}(gr E)$
which consider the summands of $c_i$ which are compositions of at most $n$ 
arrows. The $n=1$ case is given by the Yoneda product, while
when $n$ is big enough we get the whole $c_i$.
Correspondingly we have a filtration
\[0\subset H^i[1](E)\subset H^i[2](E)\subset\ldots \subset H^i(E)\]

{\bf Remark}
 The hypercohomology module of the 
complex
\[gr E\rig{\theta\wedge}gr E\otimes T_X\rig{\theta\wedge}
gr E\otimes \wedge^2T_X\rig{\theta\wedge} \ldots\]
is another interesting invariant of $E$
(compare with \cite{Simp} page 24). 
The computation in the case $E=K_X$ shows that
this should be related to the filtration above
if we twist by ${\cal O}(t)$ and sum over $t\in{\bf Z}$.

\medskip

\section{Moduli and Stability}

For simplicity we restrict in this section to the case when $X$ is an
irreducible Hermitian symmetric variety. We consider now the moduli problem of homogeneous bundles $E$
on $X$ with the same
$gr E$. Any ${\cal R}$-module $F=\oplus V_{\lambda}\otimes E_{\lambda}$
 corresponds to the dimension vector
$\alpha=(\alpha_{\lambda})\in {\bf Z}^{({\cal Q_X})_0}$
 where $\alpha_{\lambda}=\dim V_{\lambda}$.
The group $$GL(\alpha):=\prod_{{\lambda\in({\cal 
Q_X})_0}}GL(V_{\lambda})$$
 acts over $${\cal K}({\cal Q_X},\alpha):=\oplus_{a\in({\cal Q_X})_1}
Hom(V_{ta},V_{ha})$$ and over the closed subvariety
$$V_X(\alpha)\subset{\cal K}({\cal Q_X},\alpha)$$
defined by the relations in ${\cal Q}_X$. 
The affine quotient $Spec({\bf C}[V_X(\alpha)]^{GL(\alpha)})$
is a single point, represented by $F$ itself.
King (\cite{King}) considers the characters  of $GL(\alpha)$
which are given by 
$$\chi_{\sigma}(g)=\prod_{{\lambda\in({\cal 
Q_X})_0}}\det(g_{\lambda})^{\sigma_{\lambda}}$$
for $\sigma \in {\bf Z}^{({\cal Q_X})_0}$
such that $\sum_{\lambda}\sigma_{\lambda}\alpha_{\lambda}=0$.
The element $\sigma$ can also be interpreted as a homomorphism
$K_0(R$-$mod)\to {\bf Z}$ which from $E_{\lambda}$
gives $\sigma_{\lambda}$.
A function $f\in  {\bf C}[V_X(\alpha)] $
is called a relative invariant of weight $\sigma$ if
$f(g\cdot x)=\sigma(g)f(x)$, and the space of such relatively
invariant functions is denoted by 
${\bf C}[V_X(\alpha)]^{{GL(\alpha)},\sigma}$.

There is a natural character,
 that it is convenient to denote
by $\mu({\alpha})$, defined by  
$$\mu({\alpha})_{\lambda}=
c_1(F)rk(E_{\lambda})-rk(F)c_1(E_{\lambda})$$
Observe that $\mu({\alpha})(F)=
\sum_{\lambda}\alpha_{\lambda}\mu({\alpha})_{\lambda}=0$. 
For any subrepresentation $E'$ of $E\in M_X(\alpha)$ let
$gr E'=\oplus V_{\lambda}'\otimes E_{\lambda}$
with $\dim V_{\lambda}'=\alpha'_{\lambda}$,
then 
\begin{equation}
\label{muformula}
\mu(\alpha)(E')=\sum_{\lambda}\alpha'_{\lambda}\mu(\alpha)_{\lambda}=
rk E' rk F\left( \mu(F)-\mu(E')\right)
\end{equation}

Then we define 
$$M_X(\alpha):= Proj(\oplus_{n\ge 0}
{\bf C}[V_X(\alpha)]^{{GL(\alpha)},n\mu(\alpha)})$$
which is projective over $Spec({\bf C}[V_X(\alpha)]^{GL(\alpha)})$,
hence it is a projective variety. The moduli space
$M_X(\alpha)$ is the GIT quotient
 of the open set $V_X(\alpha)^{ss}$ of 
$\chi_{\mu(\alpha)}$-semistable points (\cite{King}). Different characters
give moduli spaces which are birationally equivalent to
$M_X(\alpha)$.

We collect the known results about this topic
in the following propositions. We saw that $E$ is determined by
$\theta_E\in Hom(gr E, gr E\otimes T_X)$ such that 
$\theta_E\wedge\theta_E=0$ (\thmref{starting}).

\begin{thm0}
\label{stab1} Let $E$ be a homogeneous bundle on $X$
irreducible Hermitian symmetric variety  and let
 $\alpha$
be the dimension vector corresponding to $gr E$.
 The following facts are equivalent

(i) for every $G$-invariant subbundle $K$ we have $\mu(K)\le \mu(E)$
(\underbar{equivariant semistability})

(ii) for every subbundle $K$ such that $\theta_E(gr K)\subset gr K\otimes T_X$ we have $\mu(K)\le \mu(E)$ 
(\underbar{Higgs semistability})

(iii) the 
representation $[E]$ of ${\cal Q}_X$ is $\mu(\alpha)$-semistable according
to \cite{King}, Def. 1.1. (\underbar{quiver semistability})

(iv) $E$ is a $\chi_{\mu(\alpha)}$-semistable point in $V_X(\alpha)$
for the action of $GL(\alpha)$
(\cite{King}, Def. 2.1)  (\underbar{GIT semistability})

(v) for every subsheaf $K$ we have $\mu(K)\le \mu(E)$
(\underbar{Mumford-Takemoto semistability},
 see \cite{OSS}).
\end{thm0}

{\it Proof} (i) $\Longleftrightarrow$ (ii) follows from the fact that
$F\subset E$ is $G$-invariant
iff $\theta_E(gr F)\subset gr F\otimes T_X$ .
(ii) $\Longleftrightarrow$ (iii) is straightforward from \thmref{main2},
the remark after it
and (\ref{muformula}). 
(iii)$\Longleftrightarrow$ (iv) is proved
in \cite{King}, Prop. 3.1 and Thm. 4.1. (i)$\Longleftrightarrow$ (v) is proved
in \cite{Migl} and independently in \cite{Ro} 
(this last only in the case of ${{\bf P}}^n$,
but his proof extends in a straightforward way to any $G/P$,
 see \cite{Ot}). \qedd

{\bf Remark} Migliorini shows in \cite{Migl} in the analytic
setting that conditions
(i) to (v) are equivalent to the existence of an approximate
Hermite-Einstein metric, which can be chosen invariant for
a maximal compact subgroup of $G$. He also relates the stability
to the image of the moment map.

\begin{thm0}
\label{stab2}
 Let $E$ be a homogeneous bundle on $X$
irreducible Hermitian symmetric variety and let
 $\alpha$
be the dimension vector corresponding to $gr E$.
The following facts are equivalent

(i) for every $G$-invariant proper subbundle $K$ we have 
$\mu(K)< \mu(E)$
(\underbar{equivariant} \underbar{stability})

(ii) for every proper subbundle $K$ such that
 $\theta_E(gr K)\subset gr K\otimes T_X$ we have $\mu(K)< \mu(E)$
 (\underbar{Higgs stability})

(iii) the 
representation $[E]$ of ${\cal Q}_X$ is 
$\mu(\alpha)$-stable according
to \cite{King}, Def. 1.1. (\underbar{quiver stability})

(iv) $E$ is a $\chi_{\mu(\alpha)}$-stable point in $V_X(\alpha)$
for the action of $GL(\alpha)$
(\cite{King}, Def. 1.2)  (\underbar{GIT stability})

(v) $E\simeq W\otimes E'$ where $W$ is an irreducible $G$-module and 
for every proper subsheaf $K\subset E'$ we have $\mu(K)\le \mu(E')$
(\underbar{Mumford-Takemoto stability} of $E'$,
 see \cite{OSS}).
\end{thm0}

{\it Proof} (i) $\Longleftrightarrow$ (ii) $\Longleftrightarrow$
 (iii) $\Longleftrightarrow$ (iv) are as above. (i) $\Longleftrightarrow$ (v) is proved
in \cite{Fa}.

\qedd

{\bf Remark}  The equivalence (i)$\Longleftrightarrow$ (v)
holds in the two previous theorems over any rational
homogeneous variety $X$ (for any slope $\mu_a$).

{\bf Remark} \thmref{stab1} and \thmref{stab2} extend
in a straightforward way to any $\sigma\colon K_0(R$-$mod)\to
 {\bf Z}$ such that $\sigma(gr E)=0$ at the place of $\mu(\alpha)$.

{\bf Remark} \thmref{stab2} shows that Mumford-Takemoto
stability is a stronger condition than stability in ${\cal Q}_X$.
The Euler sequence on ${\bf P}^n$ just explains this fact. Indeed
${\cal O}\otimes V$ corresponds to a stable representation
of ${\cal Q}_{{\bf P}^n}$, but it is not a Mumford-Takemoto stable bundle.
The points in $M_X(\alpha)$
parametrize  $S$-equivalent classes of semistable homogeneous bundles 
$E$ with the same $gr E$ corresponding to $\alpha$. The closed orbits
in $V_X(\alpha)^{ss}$ correspond to direct sums
$\oplus_j W_j\otimes F_j$ where $W_j$ are irreducible $G$-modules and
$F_j$ are Mumford-Takemoto stable homogeneous bundles.

When $E$ is a Mumford-Takemoto homogeneous stable bundle,
we get $W={\bf C}$ in condition (v) and
 an open set containing the corresponding point in $M_X(\alpha)$ embeds in
the corresponding Maruyama scheme of stable bundles
(see the construction of families in \S 5 of \cite{King}).
The tangent space at this point is $H^1(End E)^G$.

Observe that the irreducible bundles do not deform as homogeneous bundles
and their corresponding moduli space in the sense above is a single point.
(see \corref{irrexti}).

\begin{exa}\label{crossratio}
We describe an example of a homogeneous bundle on ${\bf P}^2$
with a continuous family of homogeneous deformations. This example 
appears already in \cite{Hi1}, ex. 1.8.7 and prop. 4.2.4.

Such example is $E=Sym^2Q(-1)\otimes{\cal S}^{2,1}V$ of rank $24$.
It is easy to compute that $H^1(End E)^{G}={\bf C}$.
The corresponding representation of the quiver associates to
\[\begin{array}{ccccc}
{\cal O}&\lef{}&Q(1)\cr
\dow{}&&\dow{}\cr
Q(-2)&\lef{}&Sym^2Q(-1)&\lef{}&Sym^3Q\cr
&&\dow{}&&\dow{}\cr
&&Sym^3Q(-3)&\lef{}&Sym^4Q(-2)\cr
\end{array}
\]
the diagram
\[\begin{array}{ccccc}
{\bf C} &\lef{}&{\bf C}\cr
\dow{}&&\dow{f_1}\cr
{\bf C}&\lef{f_4}&{\bf C}^2&\lef{f_2}&{\bf C}\cr
&&\dow{f_3}&&\dow{}\cr
&&{\bf C}&\lef{}&{\bf C}\cr
\end{array}
\]
The 4 arrows starting or ending in the middle ${\bf C}^2$
determine 4 one dimensional spaces (two kernel and two images)
which correspond to $4$ marked points in ${\bf P}^1$.
The cross-ratio of these 4 points decribes the deformation.
The generic deformation is Mumford-Takemoto stable.
If we fix the dimension vector $\alpha=(1,1,1,2,1,1,1)$
according to the diagram above, then
$$M_{{\bf P}^2}(\alpha)={\bf P}^1$$
Indeed  the character $\mu(\alpha)$ is
$72(0,-1,1,0,-2,2,0)$. We can divide by $72$ and 
the coordinate ring
$$\oplus_{n\ge 0}
{\bf C}[V_X(\alpha)]^{{GL(\alpha)},n\mu(\alpha)}) $$
is generated by 
\[S=(f_4f_1)(f_3f_2)^2\qquad\hbox{and}\qquad T=(f_4f_2)(f_3f_2)(f_3f_1)\]
(both correspond to $n=1$).
If we do not divide by $72$ then the two generators are $S^{72}$
and $T^{72}$.

There are three distinguished points.The first one
(corresponding to $S=0$) when $Im f_1=Ker f_4$. In this case there are
three different orbits where the $S$-equivalence class
contains ${\cal O}$ as direct summand.
The second one (corresponding to $S=T$) when $Im f_1=Im f_2$ or when
$Ker f_3=Ker f_4$. In this case there are
three different orbits where the $S$-equivalence class
contains $Sym^2Q(-1)$ as direct summand.
The third one (corresponding to $T=0$) when $Im f_1=Ker f_3$ or when
$Im f_2=Ker f_4$. Also in this case there are
three different orbits where the $S$-equivalence class
contains $ad V$ as direct summand.
Observe that $Im f_2=Ker f_3$ gives a nonstable situation
where  the middle row
${\bf C}\lef{f_4}Im f_2\lef{f_2}{\bf C}$
destabilizes.

There are other two particular points in $M_{{\bf P}^1}(\alpha)$
 which correspond respectively
 to $Sym^2Q(-1)\otimes{\cal S}^{2,1}V$ 
and to $ad C$ where $C$ is the rank $5$ exceptional bundle
defined by the sequence
\[0\rig{}Q(-1)\rig{}C\rig{}Sym^2Q\rig{}0\]

\end{exa}

{\bf Remark } It seems an interesting open questions
to understand when $M_X(\alpha)$ is nonempty or irreducible.

\section{Olver maps and explicit relations for Grassmannians}

The aim of this section is to make explicit in the case of Grassmannians the relations coming from $\theta\wedge\theta=0$
and the corresponding complex $H^*(gr E)$.

We restrict to the case $G=SL(V)$. 
Let $a$ (resp. $a'$, $a''$) be the Young diagram associated to
$\lambda$ (resp. $\lambda'$, $\lambda''$), so that ${\cal S}^aV$
is the representation with maximal weight $\lambda$. We have that $a'$ 
is obtained adding one box to  $a$ and we have the Pieri maps
${\cal S}^aV\otimes V\rig{}{\cal S}^{a'}V$.
These maps are defined up to a nonzero scalar multiple. Olver gave in 
the unpublished preprint \cite{Ol}
a nice description of these maps. This description was used in \cite{D},
then a proof appeared in \cite{M-O}, in the more general setting of skew Young diagrams.

It is well known that ${\cal S}^{a}V$ can be obtained as a quotient
of $Sym^{a}V:=Sym^{a_1}V\otimes\ldots\otimes Sym^{a_n}V$ (see \cite{DC-E-P} or \cite{F-H}),
namely there is the quotient map (\cite{D} 2.6)
$$\rho_{a}\colon Sym^{a}V\rig{}{\cal S}^{a}V$$
Olver's idea is to consider the Pieri maps 
at the level of   $Sym^{a}V$
and then factor through the quotient.

We follow here \cite{D}, where a different notation is used, in particular
$Sym^{\tilde a}V$ in \cite{D} is ours $Sym^{a}V$. We refer to \cite{D} for
the definition of the linear map
$\chi_{a}^{a'}\colon Sym^{a'}V\rig{} Sym^{a}V\otimes V$. This is called an 
{\it Olver map}.

\begin{thm0}\label{olverpn}{\bf (Olver, \cite{D} thm. 2.14)}
Consider the diagram
$$\begin{array}{ccc}
Sym^{a'}V&\rig{\chi_{a}^{a'}}&Sym^{a}V\otimes V\cr
\dow{\rho_{a'}}&& \dow{\rho_{a}\otimes 1}  \cr
{\cal S}^{a'}V&& {\cal S}^{a}V\otimes V\cr
\end{array}$$
Then $\chi_{a}^{a'}(\ker \rho_{a'})\subset
\ker(\rho_{a}\otimes 1)$ and $\chi_{a}^{a'}$
induces the nonzero $SL(V)$-equivariant 
$$\psi_{a}^{a'}\colon{\cal S}^{a'}V\rig{} {\cal S}^{a}V\otimes V $$
making the above diagram commutative.
\end{thm0}

 A {\it tableau} on the Young diagram $a$ is a numbering of the boxes
with the integers between $1$ and $n+1$.
A tableau is called {\it standard} if the rows are weakly increasing from the left to the right
and the columns are strictly increasing from the top to the bottom.
The {\it content} of a tableau $T$ is the function
$C_{T}\colon \{1,\ldots, n\}\to {\bf N}$ such that
$C_T(p)$ is the number of times $p$ occurs in $T$.
After a basis $e_1,\ldots,e_{n+1}$ of $V$ has been fixed, to any tableau $T$ is
associated in the natural way a tensor $T^S$ in $Sym^{a}V$
by symmetrizing the basis vectors labelled by each row. The eigenvectors for the action of the diagonal
subgroup of $SL(V)$ over ${\cal S}^{a}V$ correspond
to $\rho_{a}(T^S)$ with $T$ choosen among
the standard tableau. They form a basis of ${\cal S}^{a}V$.

\medskip

Let $K^{a}$ be the tableau obtained by filling the $i$-th row with entries equal to $i$ (it is called 
{\it canonical}
in \cite{DC-E-P});
$K^{a}$ is the only standard tableau among those with 
the same content.
The projection $\rho_{a}(K^{a S})$ is a maximal eigenvector for
${\cal S}^{a}V$ and we denote it by  $\kappa^{a}$.
Let $a'$ be obtained  from $a$ by adding a box to the 
$i$-th row
and let $a''$ be obtained from $a'$ by adding a box to the $j$-th row.
 Consider the map
$\chi_{a}^{a''}\colon Sym^{a''}V\rig{}Sym^{a}V\otimes V\otimes V$
defined as the composition
{\small  \[
Sym^{a''}V\rig{\chi_{a'}^{a''}} Sym^{a'}V\otimes V\rig{\chi_{a}^{a'}\otimes 1} 
Sym^{a}V\otimes V\otimes V
\]}

$\chi_{a}^{a''} $ also induces the nonzero $SL(V)$-equivariant morphism
$$\psi_{a}^{a''}\colon {\cal S}^{a''}V\rig{}{\cal S}^{a}V\otimes V\otimes V$$

Let $K_{i,j}^{a'}$ be the tableau on $a'$ obtained
by adding a box filled with $j$ at the $i$-th row of $K^{a}$.
We denote by $\kappa_{i,j}^{a'}$
the element $\rho_{a'}(K_{i,j}^{a S})$
\begin{prop0}
\label{olvermap} 
(i) If $i>j$ then 
$$\psi_{a'}^{a''}(\kappa^{a''})=
(a_j+1)\kappa^{a'}\otimes e_j+\sum_{h\neq j, i}\tau_h\otimes e_h$$
for some $\tau_h$.

(ii) If $i=j$ then 
$$\psi_{a'}^{a''}(\kappa^{a''})=
(a_j+2)\kappa^{a'}\otimes e_j+\sum_{h\neq j}\tau_h\otimes e_h$$
for some $\tau_h$.

(iii) If $i<j$ then 
$$\psi_{a'}^{a''}(\kappa^{a''})=\left(-\frac{(a_j+1)(a_i+1)}{a_i-a_j+j-i}\kappa_{i,j}^{a'}+\tau\right)\otimes e_i+
(a_j+1)\kappa^{a'}\otimes e_j+\sum_{h\neq i,j}\tau_h\otimes e_h$$
for some $\tau$, $\tau_h$, where  $\psi_{a}^{a'}(\tau)$ has zero coefficient in $e_j\otimes \kappa^{a}$.
\end{prop0}
{\it Proof}
In (i) and (ii) the summand $\kappa^{a'}\otimes e_j$ is obtained
with $J=(0,j)$ (see \cite{D}2.12).
In (iii) the summand   $\kappa_{i,j}^{a'}\otimes e_i$ is obtained
with $J=(0,i,j)$  while the summand
$\kappa^{a'}\otimes e_j$ is obtained with $J=(0,j)$.
\qedd

\begin{coro0}
\label{olver1} 
(i) If $i>j$ then 
$$\psi_{a}^{a''}(\kappa^{a''})=
(a_i+1)(a_j+1)\kappa^{a}\otimes e_i\otimes e_j+$$
$${\ldots\hbox{(linear combination of other basis vectors different
from $\kappa^a\otimes e_j\otimes e_i$ )}}\hfill$$

(ii) If $i=j$ then 
$$\psi_{a}^{a''}(\kappa^{a''})=
(a_j+1)(a_j+2)\kappa^{a}\otimes e_j\otimes e_j+
\ldots\hbox{(linear combination of other basis vectors)}$$

(iii) If $i<j$ then 
$$\psi_{a}^{a''}(\kappa^{a''})=
(a_i+1)(a_j+1)\kappa^{a}\otimes\left(e_i\otimes e_j-\frac{1}{a_i-a_j+j-i}e_j\otimes e_i\right)
+\ldots$$
$$\ldots\hbox{(linear combination of other basis vectors)}\hfill$$
\end{coro0}

{\bf Remark} The case (i) of \corref{olver1} does not appear if $i=j+1$ and $a_i=a_j$,
in such a case $a''$ is obtained from $a$ by adding two boxes to the same column, and the only possibility
is to add first the highest box and then the lowest one.

Now consider a bundle $E_{\lambda}=S^{\alpha}U\otimes S^{\beta}Q^*(t)$
 (as in section 5) in the Grassmannian $Gr({\bf P}^k, {\bf P}^n)$
where
$\lambda=\sum_{i=1}^n c_i\lambda_i$.  Let $p, q\in{\bf N}$.

Let 
$$n_{p,q}:=-\sum_{i=-(p-1)}^{q-1}\alpha_{k+1+i}$$
$$\lambda_{p,q}:=\lambda+n_{p,q}$$
$$\lambda_{q}:=\lambda_{1,q}=\lambda-\sum_{i=0}^{q-1}\alpha_{k+1+i}$$

We denote the corresponding morphism as
$$m_{\lambda,p,q}\colon E_{\lambda}\otimes\Omega^1\to E_{\lambda_{p,q}}$$
normalized according to \defref{defmain}.

Then $E_{\lambda_{p,q}}=S^{\alpha'}U\otimes S^{\beta'}Q^*(t)$
where $\alpha'$ is obtained from $\alpha$ by adding a box to row $p$ and
$\beta'$ is obtained from $\beta$ by adding a box to row $q$.

In the following proposition we make the relations (see \defref{rel12})
explicit for  ${\cal Q}_{Gr({\bf P}^k,{\bf P}^n)}$.
We consider $E_{\lambda''}=S^{\alpha''}U\otimes S^{\beta''}Q^*(t)$
where $\alpha''$ is obtained from $\alpha$ by adding two box to rows $p_1$, $p_2$ and
$\beta''$ is obtained from $\beta$ by adding two boxes to rows $q_1$, $q_2$. If
$p_1=p_2$ and $q_1=q_2$ then $Ext^2(E_{\lambda},E_{\lambda''})^G=0$. By the symmetry 
we may assume $p_1\le p_2$, $q_1<q_2$. Let

$$\tilde p:=\sum_{i=p_1}^{p_2-1}c_{k+1-i}+p_2-p_1=\alpha_{p_1}-\alpha_{p_2}+p_2-p_1$$

$$\tilde q:=\sum_{i=q_1}^{q_2-1}c_{k+1+i}+q_2-q_1=\beta_{q_1}-\beta_{q_2}+q_2-q_1$$

Note that ${\tilde p}=1$ if and only if $p_2=p_1+1$ and $c_{k+1-p_1}=0$. In the same way
${\tilde q}=1$ if and only if $q_2=q_1+1$ and $c_{k+1+q_1}=0$.

\begin{prop0}{\bf (Explicit relations for ${\cal Q}_{Gr({\bf P}^k,{\bf P}^n)}$)}
\label{relgras}

(i) If $p_1<p_2$, we have the subcases

(i1) ${\tilde p}\neq 1$, ${\tilde q}\neq 1$; in this case we have the two equations

{\small
$$g_{\lambda_{p_1,q_1},p_2q_2} g_{\lambda,p_1,q_1}
\left(\frac{1}{\tilde q}-\frac{1}{\tilde p}
\right)-g_{\lambda_{p_1,q_2},p_2q_1} g_{\lambda,p_1,q_2}
+g_{\lambda_{p_2,q_1},p_1q_2} g_{\lambda,p_2,q_1}=0$$}
\medskip
{\small
$$g_{\lambda_{p_1,q_1},p_2q_2} g_{\lambda,p_1,q_1}
\left(\frac{1}{{\tilde p}{\tilde q}}-1
 \right)+$$
$$-g_{\lambda_{p_1,q_2},p_2q_1} g_{\lambda,p_1,q_2}
\left(
\frac{1}{\tilde p}
\right)
-g_{\lambda_{p_2,q_1},p_1q_2} g_{\lambda,p_2,q_1}
\left( \frac{1}{\tilde q} \right)
+g_{\lambda_{p_2,q_2},p_1q_1} g_{\lambda,p_2,q_2}=0$$ }

(i2) ${\tilde p}=1$ and ${\tilde q}\neq 1$; in this case   $\lambda_{p_2,q_1}, \lambda_{p_2,q_2}$
do not exist and we have the single equation
$$g_{\lambda_{p_1,q_1},p_2q_2} g_{\lambda,p_1,q_1}\left(\frac{1}{{\tilde q}}-1
 \right)-g_{\lambda_{p_1,q_2},p_2q_1} g_{\lambda,p_1,q_2}=0$$

(i3) ${\tilde p}\neq 1$ and ${\tilde q}=1$; in this case   $\lambda_{p_1,q_2}, \lambda_{p_2,q_2}$
do not exist and we have the single equation
$$g_{\lambda_{p_1,q_1},p_2q_2} g_{\lambda,p_1,q_1}\left(1-\frac{1}{{\tilde p}}
 \right)+g_{\lambda_{p_2,q_1},p_1q_2} g_{\lambda,p_2,q_1}=0$$

(i4) ${\tilde p}={\tilde q}=1$; in this case  only $\lambda_{p_1,q_1}$
survives and there are no equations at all.  Hille counterexample (see \exref{g13})
fits this case.

\medskip

(ii) If $p_1=p_2$ we have the subcases 

(ii1) ${\tilde q}\neq 1$; in this case  we have the equation

$$g_{\lambda_{p_1,q_1},p_1q_2} g_{\lambda,p_1,q_1}\left(\frac{1+\tilde q}{\tilde q}\right)-
g_{\lambda_{p_1,q_2},p_1q_1} g_{\lambda,p_1,q_2}=0$$
               
(ii2) ${\tilde q}=1$; in this case  we have the equation

$$g_{\lambda_{p_1,q_1},p_1q_2} g_{\lambda,p_1,q_1} =0$$

 \end{prop0}
 
 {\bf Proof}
Let $p_1<p_2$. Consider that

$$m_{\lambda_{p_1,q_1},p_2q_2}\wedge m_{\lambda,p_1,q_1}(n_{p_1q_2}\wedge n_{p_2q_1}\otimes v_{\lambda})=
m_{\lambda_{p_1,q_1},p_2q_2} \left(n_{p_1q_2}\otimes m_{\lambda,p_1,q_1}(n_{p_2q_1}\otimes v_{\lambda})\right)+$$
$$-
m_{\lambda_{p_1,q_1},p_2q_2} \left(n_{p_2q_1}\otimes m_{\lambda,p_1,q_1}(n_{p_1q_2}\otimes v_{\lambda})\right)
= \left(-\frac{1}{\tilde p}+\frac{1}{\tilde q}\right)v_{\lambda''}$$
(the last equality by \corref{olver1}).
In the same way  if ${\tilde q}\neq 1$
$$m_{\lambda_{p_1,q_2},p_2q_1}\wedge m_{\lambda,p_1,q_2}(n_{p_1q_2}\wedge n_{p_2q_1}\otimes v_{\lambda})=
-v_{\lambda''}$$

Moreover  if ${\tilde p}\neq 1$
$$m_{\lambda_{p_2,q_1},p_1q_2}\wedge m_{\lambda,p_2,q_1}(n_{p_1q_2}\wedge n_{p_2q_1}\otimes v_{\lambda})=
v_{\lambda''}$$
Besides
$$m_{\lambda_{p_2,q_2},p_1q_1}\wedge m_{\lambda,p_2,q_2}(n_{p_1q_2}\wedge n_{p_2q_1}\otimes v_{\lambda})=0$$
Now by computing the left side of the relation (\ref{relationmg}) on $n_{p_1q_2}\wedge n_{p_2q_1}\otimes v_{\lambda}$
we get the first equation of (i1).

In the same way,  computing the left side of the relation (\ref{relationmg}) on 
$n_{p_1q_1}\wedge n_{p_2q_2}\otimes v_{\lambda}$
we get the second equation of (i1). The other subcases of (i) are particular cases of (i1).
(ii) is analogous.\qedd

\medskip

 {\bf Remark} The number of equations obtained in \propref{relgras} measures exactly 
the dimension of $Ext^2(E_{\lambda}, E_{\lambda''})^G$, which can be $2$, $1$, or $0$. An interesting consequence of
\propref{relgras} is that (with the assumptions in (i)) there is no indecomposable homogeneous bundle
on $Gr({\bf P}^k,{\bf P}^n)$ such that  its quiver representation has support equal to the  parallelogram
with vertices
$E_{\lambda}$, $E_{\lambda_{p_1,q_2}}$, $E_{\lambda_{p_2,q_1}}$, $E_{\lambda''}$.
The first consequence is that on  $Gr({\bf P}^1,{\bf P}^3)$ every homogeneous
bundle $E$ such that $gr E=\Omega^1\oplus\Omega^2\oplus\Omega^3$ decomposes.   On the other hand
an indecomposable homogeneous bundle
such that  its quiver representation has support equal to the  parallelogram
with vertices
$E_{\lambda}$, $E_{\lambda_{p_1,q_1}}$, $E_{\lambda_{p_2,q_2}}$, $E_{\lambda''}$
exists if and only if ${\tilde p}={\tilde q}$.
The first nontrivial example is, on the Grassmannian of lines in ${\bf P}^3={\bf P}(V)$ ,
the cohomology bundle $E$ of the monad
$${\cal O}(-2)\rig{}{\cal S}^{2,2}V\rig{}{\cal O}(2)$$
which has $gr E={\cal O}\oplus\Omega^1\oplus\Omega^1(2)\oplus (Sym^2U\otimes Sym^2Q)$.

\medskip
 
 \begin{coro0}{\bf (Explicit relations for ${\cal Q}_{{\bf P}^n}$)}
 \label{commutpn}
 In the case of ${\bf P}^n$
 the category of homogeneous bundles is equivalent to the category of
 representations of ${\cal Q}_{{\bf P}^n}$ with the commutativity relations.
 \end{coro0} 

{\bf Proof}
Put $p_1=p_2=1$ in \propref{relgras} and get 

$$g_{\lambda_{q_1},q_2} g_{\lambda,q_1}
\left(\frac{1+\tilde q}{\tilde q}\right)-
g_{\lambda_{q_2},q_1} g_{\lambda,q_2}=0$$  

unless ${\tilde q}=1$. 

Denoting
$$h_{\lambda,i}:=(c_i+1)(c_{i-1}+c_i+2)\cdots (c_2+\ldots +c_i+i-1)f_{\lambda,i}$$
we get a functor from the quiver ${\cal Q}_{{\bf P}^n}$ with the relations that we have defined
to the same quiver but with  the commutativity relations 
 $$h_{\lambda_{q_1},q_2} h_{\lambda,q_1}-h_{\lambda_{q_2},q_1} h_{\lambda,q_2}=0$$
This functor gives the desired equivalence.
\qedd

\bigskip

{\bf Address (of both authors)}: Dipartimento di
Matematica "U.Dini", Viale Morgagni 67/A,  50134 Firenze,
Italy

\smallskip

{\bf E-mail addresses}: ottavian@math.unifi.it,
rubei@math.unifi.it

\end{document}